\documentclass[12pt, a4paper, oneside]{article}
\usepackage[body={16cm, 24cm},left=2cm,right=2cm]{geometry}
\usepackage{float}

\usepackage{natbib}

\usepackage{amsmath,amssymb,amsthm,latexsym}

\usepackage{graphicx}
\usepackage{setspace}	% Line spacing
\usepackage{pifont}
\usepackage{color}
\usepackage{fancyhdr}
\usepackage{hyperref}
\newtheorem{thm}{Theorem}
\newtheorem{lem}[thm] {Lemma}
\newtheorem{rem}[thm]{Remark}
\newtheorem{cor}[thm]{Corollary}
\newtheorem{prop}[thm]{Proposition}

\newcommand{\keywords}[1]{{\scriptsize \noindent \textbf{KEY WORDS AND PHRASES:}} {#1}\\}

\newcommand{\id} {\ensuremath{\displaystyle{\mathop {=} ^d}}}

\newcommand{\field}[1]{\mathbb{#1}}
\newcommand{\real}{\ensuremath{{\field{R}}}}
\newcommand{\mc}[1]{{\ensuremath{\mathcal{#1}}}}

\newcommand{\sumab}[2]{\ensuremath{\sum\limits_{#1}^{#2}}}
\newcommand{\intab}[2]{\ensuremath{\int_{#1}^{#2}}}

\newcommand{\intunit}{\ensuremath{\int_{0}^{1}}}

\newcommand{\limit}[1]{\ensuremath{\displaystyle {\lim_{#1 \rightarrow{\infty}}}}}

\newcommand{\conv}[1]{\ensuremath{\, \displaystyle {\mathop {\longrightarrow}_{n \rightarrow \infty} ^{#1}}}\, }

\newcommand{\kdivn}{\frac{k}{n}}
\newcommand{\ndivk}{\frac{n}{k}}

\DeclareMathOperator*{\argmax}{arg\,max}

\title{{\LARGE A general estimator for the right endpoint}\\
{\large with an application to supercentenarian women's records}\thanks{Funded by FCT - Funda\c{c}\~{a}o para a Ci\^{e}ncia e a Tecnologia, Portugal, through the project UID/MAT/00006/2013}
}
\author{Isabel Fraga Alves (\emph{\small CEAUL, University of Lisbon})
\and Cl\'{a}udia Neves (\emph{\small University of Reading, UK})
 \and Pedro Ros\'{a}rio (\emph{\small CEAUL, University of Lisbon})
}

\date{ }
\begin{document}
\maketitle

%\begin{abstract}
\abstract{We extend the setting of the right endpoint estimator introduced in Fraga Alves and Neves (Statist. Sinica 24:1811--1835, 2014) to the broader class of light-tailed distributions with finite endpoint, belonging to some domain of attraction induced by the extreme value theorem. This stretch enables a general estimator for the finite  endpoint, which does not require estimation of the (supposedly non-positive) extreme value index. A new testing procedure for selecting max-domains of attraction also arises in connection with   asymptotic properties of the general endpoint estimator. The  simulation study conveys that the general endpoint estimator is a valuable complement to the most usual endpoint estimators, particularly when the true extreme value index stays above $-1/2$, embracing the most common cases in practical applications. An illustration is provided  via an extreme value analysis of supercentenarian women data.}
%\\\\

\keywords{ Extreme value theory,  Semi-parametric estimation, Tail estimation,  Regular variation,  Monte Carlo simulation, Human lifespan}
% \PACS{}
 %\subclass{62G32 \and 62F10}
%\end{abstract}

%==================================================================
\section{Introduction}
\label{SecIntro}

The extreme value theorem  \citep[with contributions from][]{FisherTippett:28, G:43, deHaan:70} and its counterpart for exceedances above a threshold \citep{BdH:74} ascertain that inference about rare events can be   drawn on the larger (or lower) observations in the sample. While restricting attention to the large rare events, the theoretical framework provided by the extreme value theorem reads as follows. If a non-degenerate limit  $G$ is achieved by the distribution function (d.f.) of the partial maxima $X_{n,n}$
of a sequence $\{X_n\}_{n\geq 1}$ of independent and identically distributed (i.i.d.) random
variables (r.v.) with common d.f. $F$, and  if  there exist $a_n>0$ and $b_n\in \real $ such that
$
	\lim_{n\rightarrow\infty}F^n(a_n\,x+b_n)= G(x),
$
for every continuity point of $G$, then $G$ is one of the  three distributions
\begin{eqnarray}
\label{Gumb} \Lambda(x) &=& \exp\{-\exp(-x)\},\quad x\in\real,\\
\label{Frech}	\Phi_\alpha(x)&=&\exp\{-x^{-\alpha}\},\;x>0,\quad \alpha>0,\\
\nonumber \Psi_\alpha(x)&=&\exp\{-(-x)^\alpha\},\;x<0,\quad \alpha>0.
\end{eqnarray}
These can be nested in the Generalized Extreme Value (GEV) d.f.
\begin{equation}\label{GEV}
G_{\gamma}(x):= \exp\{-(1+\gamma x)^{-1/\gamma}\}, \; 1+\gamma x>0,\;
\gamma \in \mathbb{R}.
\end{equation}
We then say that $F$ is in the (max-)domain of attraction of $G_\gamma$,  for some extreme value index (EVI) $\gamma \in \real$ [notation: $F \in \mathcal{D}(G_{\gamma}) $]. For $\gamma=0$, the right-hand side of \eqref{GEV} is read as $\exp\left(-e^{-x}\right)$.
The theory of regular variation \citep{Binghametal:87,deHaan:70, deHaanFerreira:06}, provides necessary and sufficient conditions for $F\in \mc{D}(G_{\gamma})$. Let $U$ be the tail quantile function defined by the generalized inverse of $1/(1-F)$, i.e.
$
U(t):=   F^{\leftarrow} \bigl( 1-1/t\bigr)$, for $ t>1.
$
Then, $F\in \mc{D}(G_{\gamma})$ if and only if there exists a positive  measurable function $a(\cdot)$ such that the condition of \emph{extended regular variation}
\begin{equation}\label{ERVU}
	\limit{t}\,\frac{U(tx)-U(t)}{a(t)}= \frac{x^{\gamma}-1}{\gamma},
\end{equation}
holds for all $x>0$ [notation: $U\in ERV_{\gamma}$].
The limit  in  \eqref{ERVU} coincides with the $U$-function of the Generalized Pareto distribution (GPD), with distribution function $1+\log G_\gamma$. Hence,   for extrapolating beyond the range of the available observations, the statistics of extremes will be exclusively focused on those observations over a sufficiently high threshold. Then the excesses above this threshold are expected to behave as observations drawn from the GPD.

 The right endpoint of the underlying distribution function $F$ is defined as
\begin{equation*}
	x^F:= \sup\{x:\, F(x)<1\} \leq \infty,
\end{equation*}
which in terms of high quantiles  is given by  $x^F= \lim_{t\rightarrow \infty} U(t)= U(\infty)$.
For estimating the right endpoint $x^F$ we will follow a semi-parametric approach, that  is, our focus is on the domain of attraction rather than on the limiting GEV distribution.
We will also assume that $k$ is  an \emph{intermediate} sequence of positive integers $k=k_n$ such that $k \rightarrow \infty$ and $k/n \rightarrow 0$, as $n\rightarrow \infty$. This is our large sample assumption for the moment. Other mild yet reasonable conditions in the context of extreme value estimation will come forth in section \ref{SecEstim}, which essentially convey suitable bounds on the intermediate sequence $k_n$.

This paper deals with a unifying semi-parametric  approach to the problem of estimating the finite right endpoint $x^F$ when $F$ belongs to some domain of attraction where a finite endpoint is admissible, more formally $F\in \mathcal{D}(G_{\gamma})_{\gamma\leq 0}$.  We  term this estimator $\hat x^F$  the general endpoint estimator. We will provide evidence that despite $\hat x^F$ not being asymptotically normal for all values of $\gamma <0$ (a drawback if one wishes to construct confidence intervals) it proves nonetheless to be a valuable tool in terms of applications. One of the most obvious estimators of the right endpoint is the sample maximum. In fact, \citet{deHaanFerreira:06} point out in their Remark 4.5.5 that using the sample maximum $X_{n,n}$ to estimate $x^F$ in case $\gamma<-1/2$ is approximately equivalent to using the moment related estimator for the  endpoint. The striking feature of the general endpoint estimator is that it avoids the nuisance of changing ``tail estimation-goggles'' each time we are dealing with yet another sample, possibly from a distribution in a different domain of attraction. We exemplify this point by referring the study by \citet{EM:08}, which could well benefit from using the same endpoint estimator at all instances, in all athletic events. This freedom of constraint about $\gamma\leq 0$ motivates the present general estimator, alongside with its preceding application to the long jump records in \citep{FAHN:13}.

The outline of the paper is as follows. Section \ref{SecMethods} introduces the general estimator  $\hat x^F$ and its theoretical assumptions, aligned with the usual semi-parametric framework. Large sample results for  $\hat x^F$  are presented in Section \ref{SecEstim}, as well as a new test statistic based on  $\hat x^F$  aimed at selecting max-domains of attraction.
Section \ref{SecSims} is dedicated to a comparative study via simulation, involving common parametric and semi-parametric inference approaches in extremes.
 Section \ref{CaseStudy} provides an illustration of the exact behaviour of the general endpoint estimator, using the \emph{supercentenarian women} data set. Here we consider two alternative settings: estimation of the right endpoint with a link to the EVI estimation, estimation of the endpoint when this link to the EVI is broken and finally in Section \ref{SecConclusion}, we list several concluding remarks. Appendix \ref{SecProofs} encloses all the proofs of the large sample results in Section \ref{SecEstim} and  Appendix  \ref{RB2_MC} encompasses the  finite sample properties for a consistent reduced bias estimator  of the endpoint, for an EVI in $(-1/2,0)$, being compared with POTML and Moment methodologies.

\section{Semi-parametric approach to endpoint estimation}
\label{SecMethods}

We now introduce some notation. Let $F$ be the d.f. of the r.v. $X$ and $X_{1,n}\leq X_{2,n}\leq \ldots \leq X_{n,n}$ be the $n$-th ascending order statistics (o.s.) associated with the sample $X_1,\ldots,X_n$ of $n$ i.i.d. copies of $X$. We assume $F\in \mathcal{D}(G_{\gamma})$, for some  $\gamma \leq 0$, and that $x^F<\infty$.

Several estimators for the right endpoint $x^F$ of a light-tailed distribution attached to
an EVI $\gamma <0$ are available in the literature  \citep[e.g.][]{Hall:82, Caietal:13, deHaanFerreira:06}. These estimators often bear on the extreme value condition \eqref{ERVU} with $x=x(t)\rightarrow \infty$, as $t\rightarrow \infty$: since $\gamma <0$ entails that $\lim_{t\rightarrow \infty}U(t)= U(\infty)$ exists finite, then relation \eqref{ERVU} rephrases as
\begin{equation*}
\limit{t} \frac{U(\infty)-U(t)}{a(t)}= -\frac{1}{\gamma}.
\end{equation*}
A valid estimator for the right endpoint $x^F= U(\infty)$ thus arises by making $t=n/k$ in the approximate equality
$
U(\infty)\approx U(t)-a(t)/\gamma\,,
$
replacing $U(n/k)$, $a(n/k)$ and $\gamma$ by suitable consistent estimators, i.e.
\begin{equation*}
\hat{x}^*= \hat{U}\bigl(\ndivk \bigr)-\frac{\hat{a}\bigl(\ndivk \bigr)}{\hat{\gamma}}
\end{equation*}
\citep[cf. Section 4.5 of][]{deHaanFerreira:06}. Typically we consider the class of endpoint estimators
\begin{equation}\label{EPClass}
\hat{x}^*=  X_{n-k,n}-\frac{\hat{a}\bigl(\ndivk \bigr)}{\hat{\gamma}}.
\end{equation}
There is however one estimator for the right endpoint $x^F$ that does not depend on the estimation of the EVI $\gamma$. This estimator, introduced in \citet{FAN:14}, is primarily tailored for distributions with finite right endpoint in the Gumbel domain of attraction. The study of consistency and asymptotic distribution of this same endpoint estimator is the main objective in this paper, while it aims at covering the whole scenario in extremes, thus providing a unified estimation procedure for the right endpoint in the case of $\gamma\leq 0$.

The general right endpoint estimator from \citet{FAN:14} is  defined as
 \begin{equation}\label{EPEst}
	\hat{x}^F:=X_{n,n}+X_{n-k,n}- \frac{1}{\log 2}\sumab{i=0}{k-1}\log \Bigl(1+\frac{1}{k+i}\Bigr) X_{n-k-i,n}\,\,.
\end{equation}
With $a_{i,k}:= \log\left(\frac{k+i+1}{k+i}\right)/\log 2$, the endpoint estimator $\hat{x}^F$ in \eqref{EPEst} can be expressed in the equivalent form
  \begin{equation}\label{EPEst_alter}
	\hat{x}^F:= X_{n,n}+ \sum_{i=0}^{k-1}a_{i,k}\, \left(X_{n-k,n}-X_{n-k-i,n}\right) \quad \mbox{ with }
 \quad\sum_{i=0}^{k-1} a_{i,k}=1.
\end{equation}
From the non-negativeness of the weighted spacings in the sum in \eqref{EPEst_alter}, it is clear that estimator $\hat{x}^F$ is greater than $X_{n,n}$, which constitutes a major advantage to the usual semi-parametric right endpoint estimators  in  the Weibull max-domain of attraction. Therefore, the estimator $\hat{x}^F$ defined in \eqref{EPEst} can be seen as a real asset in the context of semi-parametric estimation of the finite right endpoint, embracing all distributions connected with a non-positive EVI $\gamma$, which gains by far a broader spectrum of application to the usual alternatives.

\section{Endpoint estimation and testing} %==============Sec Main =====================================
\label{SecEstim}

This section contains the  main results of the paper, giving accounts of strong consistency and sometimes asymptotic normality (we will see that the limiting normal distribution is only attained if $\gamma <-1/2$) of the general endpoint estimator $\hat{x}^F$ defined in \eqref{EPEst}. A second order reduced bias version of $\hat{x}^F$ is also devised. Additionally, we provide the asymptotic framework for a statistical test aimed at discriminating between max-domains of attraction. The new test statistic  builds on the general endpoint estimator $\hat{x}^F$. All the proofs are postponed to the Appendix \ref{SecProofs}.

\begin{prop}\label{PropCons}%--------------------------Proposition--------------------------
	Suppose $x^F$ exists finite. Assume that the extended regular variation property \eqref{ERVU} holds with $\gamma\leq 0$. If $k=k_n\rightarrow \infty$, $k_n/n \rightarrow 0$, as $n\rightarrow \infty$, then the following almost sure convergence holds with respect to $\hat{x}^F$ defined in \eqref{EPEst}:
	\begin{equation*}
		\hat{x}^F \conv{a.s.} x^F,
	\end{equation*}
	then $\hat{x}^F$ is a consistent estimator for $x^F<\infty$, i.e. $\hat{x}^F \conv{p} x^F$.\\
\end{prop}
Note that if $F\in \mc{D}(G_{\gamma})$ with $\gamma>0$, then $\hat{x}^F$ converges almost surely to infinity.
We now require a second order refinement of condition \eqref{ERVU} and auxiliary second order conditions  in order to have a grasp at the speed of convergence in \eqref{ERVU}. In particular, we assume there exists a positive or negative
  function $A_0$ with $\lim_{t \rightarrow \infty}A_0(t)=0$ such that
for each $x>0$,
  \begin{equation}\label{2ERVU}
    \limit{t}\frac{\frac{U(tx)-U(t)}{a_0(t)}-\frac{x^{\gamma}-1}{\gamma}}{A_0(t)}= \Psi^{\star}_{\gamma,\rho}(x),
\end{equation}
where $\rho$ is a non-positive parameter and with
\begin{equation*}
\Psi^{\star}_{\gamma,\rho}(x):= \left\{
                              \begin{array}{ll}
                                \frac{x^{\gamma+\rho}-1}{\gamma+\rho}, & \mbox{ } \gamma+\rho \neq 0, \, \rho <0, \\
                                \log x, & \mbox{ }  \, \gamma+\rho= 0, \, \rho <0,\\
                                \frac{1}{\gamma}\,x^{\gamma} \log x, & \mbox{ }\rho =0\neq \gamma ,\\
                                \frac{1}{2}\, (\log x)^2, & \mbox{ } \gamma= \rho=0,
                              \end{array}
                            \right.
\end{equation*}
\begin{equation*}
	a_0(t):= \left\{
                              \begin{array}{ll}
                                a(t)\bigl(1-A_0(t)\bigr), & \mbox{ } \rho <0, \\
                                a(t)\bigl(1-A_0(t)/\gamma \bigr), & \mbox{ }  \, \rho= 0\neq \gamma,\\
                                a(t), & \mbox{ } \gamma= \rho=0.
                              \end{array}
                            \right.
\end{equation*}
Moreover, $|A_0|\in RV_{\rho}$ and
\begin{equation}\label{2RVa}
\limit{t} \frac{\frac{a_0(tx)}{a_0(t)}-x^{\gamma}}{A_0(t)}= x^{\gamma}\,\frac{x^{\rho}-1}{\rho},
\end{equation}
for all $x>0$  \citep[cf. Theorem 2.3.3 and Corollary 2.3.5  of][]{deHaanFerreira:06}. Denote $U(\infty):= \lim_{t\rightarrow \infty}U(t) \,\, (=x^F)$; if \eqref{2ERVU} holds with $\gamma<0$ then,  provided $x=x(t)\rightarrow \infty$,
\begin{equation}\label{2ERVinfty}
	\limit{t}\frac{\frac{U(\infty)-U(t)}{a_0(t)}+\frac{1}{\gamma}}{A_0(t)}= \Psi^{\star}_{\gamma,\rho}(\infty):= -\frac{1}{\gamma+ \rho}\,I_{\{\rho <0\}}
	\end{equation}
by similar arguments of Lemma 4.5.4 of \cite{deHaanFerreira:06}, with $I_A$ denoting the indicator function which is equal to $1$ if $A$ holds true and is equal to zero otherwise.

%%%%%%%%%%%%%%%%%%%%%%
\begin{thm}\label{ThmMain}
Let $F$ be a d.f. in the Weibull domain of attraction, i.e., $F\in \mathcal{D}(G_\gamma)$ with $\gamma <0$. Suppose $U$ satisfies condition \eqref{2ERVU} with $\gamma <0$ and, in this sequence,  assume that \eqref{2ERVinfty} holds. We define
\begin{equation}\label{aga}
	h(\gamma):=\frac{1}{\gamma}\Bigl(\frac{2^{-\gamma}-1}{\gamma\, \log 2}+ 1\Bigr).
\end{equation}
If the intermediate sequence $k=k_n$ is such that $\sqrt{k_n}\,A_0(n/k_n) \rightarrow \lambda^* \in \real$, then
\begin{equation*}
	k^{\min(-\gamma,1/2)}\, \Bigl(\frac{\hat{x}^F-x^F}{a_0\bigl(\ndivk\bigr)}-h(\gamma)\Bigr)
	\conv{d} W \, I_{\{\gamma\geq -1/2\}}+\Bigl(N- \lambda^*\, b_{\gamma, \rho}\Bigr)\, I_{\{\gamma \leq -1/2\}},
\end{equation*}
where $W$ is a max-stable Weibull r.v.,  with d.f. $\exp\{ -(\gamma x)^{-1/\gamma} \}$ for $x<0 $, $N$ is a normal r.v. with zero mean and variance given by
\begin{equation}\label{VarN}
	Var(N)= 1+\frac{2}{\gamma \,(\log 2)^2}\Bigl(\frac{2^{-(2\gamma+1)}-1}{2\gamma+1}-\frac{2^{-(\gamma+1)}-1 }{\gamma+1}+ \frac{\log 2}{\sqrt{2}}\, (2^{-\gamma}-1)\Bigr).
\end{equation}
and $b_{\gamma, \rho}$ is defined as
\begin{equation*}
	b_{\gamma, \rho}:= \frac{1}{\log 2}\intab{1/2}{1}\Psi^{\star}_{\gamma,\rho}\bigl(\frac{1}{2s}\bigr)\,\frac{ds}{s}=\left\{
                              \begin{array}{ll}
                                \frac{1}{\gamma+\rho}\Bigl(\frac{1}{\log 2}\frac{1-2^{-(\gamma+\rho)}}{\gamma+\rho}-1\Bigr) , & \mbox{ }  \rho <0, \\
                                \frac{1}{\gamma^3\log 2}\Bigl(2^{-\gamma}(\log \,2^\gamma+1)-1\Bigr), & \mbox{ }  \rho= 0.
			 \end{array}
                            \right.
\end{equation*}
Moreover, the r.v.s  $W$ and $N$ are independent.
\end{thm}
\begin{rem}
The same normalization by $(a_0(n/k))^{-1}$, with respect to $\gamma=0$, is obtained in \citet{FAN:14} towards the Gumbel limit.
\end{rem}

\begin{cor}\label{MainCor}
Under the conditions of Theorem \ref{ThmMain},
\begin{equation*}
\frac{\sqrt{k}\,\Bigl(\frac{\hat{x}^F-x^F}{a_0\bigl(\ndivk\bigr)}-h(\gamma)\Bigr)}{k^{(\gamma+1/2)^+}}\, \conv{d} R,
\end{equation*}
where $a^+:= \max (a,0)$ and $R$ is a random variable with the following characterization:
\begin{enumerate}
\item Case $-1/2<\gamma <0$: $R$ is max-stable Weibull, with d.f. $\exp\{ -(\gamma x)^{-1/\gamma} \}$ for $x<0 $, with mean $\Gamma(1-\gamma)/\gamma$ and variance equal to $\gamma^{-2}\bigl(\Gamma(1-2\gamma)- \Gamma^2(1-\gamma) \bigr)$. Here and throughout, $\Gamma(.)$ denotes the gamma function, i.e. $\Gamma(a)= \int_0^\infty t^{a-1}e^{-t}\,dt$, $a>0$.
\item Case $\gamma <-1/2$: $R$ has normal distribution with mean $-\lambda^* b_{\gamma,\rho}$ and variance given in \eqref{VarN}.
\item Case $\gamma=-1/2$:  $R$ is the sum of the two cases above, taken as independent components, which yields a random part with mean
   $ \Gamma(1/2)-\lambda^* b_{-1/2,\rho}=\sqrt{\pi}- \lambda^* b_{-1/2,\rho}$
    and variance $
    5-\pi + 4 \bigl[ 1+ (1/\sqrt{2}-1)(2+\log 2)/\log 2 \bigr] / \log 2$.
\end{enumerate}
\end{cor}

\begin{rem}\label{RB}
The function $h(\gamma)$ is monotone decreasing for all $\gamma<0$.
Taking into account the statement of Theorem \ref{ThmMain}, an adaptive reduced bias estimator based on the general estimator $\hat x_F$ is given by
$ \hat x^F_{RB1}=\hat x_F\, - \,h(\hat\gamma)\hat a (n/k)$,  with  consistent estimators $\hat\gamma$ and $\hat a (n/k)$. The dominant component of the bias comes from the scale function $a(n/k)$ which, in case $\gamma$ is close to $0$, determines  a very slow convergence. We have conducted several simulations in this respect, indicating that this bias correction (of first order) has a very limited effect.
 \end{rem}

In addition, we consider an adaptive second order reduced bias estimator developed on the general estimator $\hat x_F$ and supported on the asymptotic statement in Theorem \ref{ThmMain} for $\gamma \in (-1/2,0)$. The limiting Weibull random variable $W$ has a non-null mean equal to $\Gamma (1-\gamma)/ \gamma$. We note  that, for small negative values of $\gamma$, the convergence of the normalized general estimator $\hat x_F$ towards the Weibull limit can be very slow since it is essentially governed by the function $a \in RV_\gamma$ and by the power transform $k^\gamma$.  Formally, the general endpoint estimator $\hat x^F$ satisfies the distributional representation
\begin{equation*}
\hat x^F = x^F + h(\gamma) a_0\bigl(\ndivk\bigr)\,+\,  a_0\bigl(\ndivk\bigr) k^\gamma W + o_p\Bigl( a_0\bigl(\ndivk\bigr) k^\gamma \Bigr)\,,
\end{equation*}
with $h(\gamma)$ defined in \eqref{aga}. We thus develop an adaptive second order reduced bias estimator as follows:
\begin{equation}\label{RB2}
  \hat x^F_{RB2}=\hat x^F_{RB1} \,-\, \frac{\Gamma(1-\hat\gamma)}{\hat\gamma}\,\hat a_0(\frac{n}{k})\, k^{\hat\gamma}
\end{equation}
(see Remark \ref{RB} for the definition of $\hat{x}^F_{RB1}$). Furthermore, an approximated $100(1-\alpha)\%$-confidence upper bound for $x^F$ is given by
\begin{equation}\label{up_bound}
x^F \,<\, \hat x^F-\hat a_0(\frac{n}{k}) \bigl[ h(\hat\gamma)+  k^{\hat\gamma} \,\, q_{\alpha}\bigr] \,,
\end{equation}
with estimated $\alpha$-quantile  of the Weibull limit distribution $q_{\alpha}:= (-\log  \alpha)^{-\hat \gamma}/\hat \gamma$.\\

\noindent
In practice, it is often advisable to perform statistical tests on the EVI sign so as to  prevent against  an actual infinite endpoint. In Section \ref{CaseStudy}, the testing procedures by \citet{NevesPFA:06} and \citet{NevesFA:07} are applied with independent interest from the particular EVI estimation problem inherent to the endpoint estimation. The hypothesis-testing problem regarding the suggested max-domain of attraction selection is stated as follows:
\begin{equation}\label{TestNeq0}
	H_0:\, F\in \mc{D}(G_0) \quad vs \quad H_1:\, F\in \mc{D}(G_\gamma)_{\gamma \neq 0}.
\end{equation}
We will introduce another test statistic for tackling this testing problem. The new statistic $G_{n,k}$ arises in connection with the general endpoint estimator  $\hat x^F$, thus in complete detachment of any extreme value index estimation procedure. It is given by
\begin{equation}\label{G}
G_{n,k}:={  \hat x^F- X_{n-k,n}   \over  X_{n-k,n}-X_{n-2k,n}   }.
\end{equation}
The next Theorem comprises the testing rule and ascertains consistency of the new testing procedure with prescribed significance level $\alpha$.

%%%%%%%%%%%%%%%%%%%%%%
\begin{thm}\label{Thm.Test}
Assume $F \in\mathcal{D}(G_\gamma)$, for some $\gamma \in \real$. Furthermore assume that the tail quantile function $U$ satisfies the second order conditions \eqref{2ERVU} up to \eqref{2RVa}. We define
\begin{equation}\label{Gnkg}
G_{n,k}^*(0):=\log 2\, G_{n,k}-\bigl( \log k + \frac{\log 2}{2}\bigr),
\end{equation}
where $G_{n,k}$ is given in \eqref{G}. If  $k=k_n$ is an intermediate sequence such that $\sqrt{k}\, A_0(n/k) \rightarrow \lambda^* \in \real$, as $n\rightarrow \infty$, then
\begin{itemize}
\item $G_{n,k}^*(0)\conv{d} Z$, if $\gamma =0$. Here  $Z$ has Gumbel distribution function $\Lambda = G_0$;
\item $G_{n,k}^*(0)\conv{P} +\infty$ , if $\gamma >0$;
\item $G_{n,k}^*(0)\conv{P} -\infty$ , if $\gamma <0$.
\end{itemize}
\end{thm}

Denoting by $\xi_{p}:=-\log(-\log(p))$ the $p$-quantile of the Gumbel distribution, a critical region for the two-sided test postulated in \eqref{TestNeq0}, at an approximate $\alpha$-level, is deemed by Theorem \ref{Thm.Test}. The statement is that we reject $H_0$ if  either $G_{n,k}^*(0)\leq \xi_{\alpha/2}$ or  $ G_{n,k}^*(0)\geq \xi_{1-\alpha/2}$.
Theorem \ref{Thm.Test} also allows testing for the one-side counterparts:
\begin{itemize}
  \item $H_0: F\in{\mathcal D}(G_0)\quad vs \quad H_1: F\in{\mathcal D}(G_\gamma)_{\gamma>0}$,\\
      thus rejecting $H_0$ in favour of heavy-tailed distributions, if $G_{n,k}^*(0)\geq \xi_{1-\alpha}$;
  \item $H_0: F\in{\mathcal D}(G_0)\quad vs \quad H_1: F\in{\mathcal D}(G_\gamma)_{\gamma<0}$,\\
      by rejecting $H_0$ in favour of short-tailed distributions, if $G_{n,k}^*(0)\leq \xi_{\alpha}$.
\end{itemize}

Denoting the power function for the testing problem \eqref{TestNeq0} (and subsequent one-sided alternatives) by
$\beta_n(\gamma):=P_\gamma[\mbox{reject } H_0]$, it follows immediately from Theorem  \ref{Thm.Test} that all the designed tests are consistent tests since, as $n\rightarrow \infty$, $\beta_n(\gamma)|_{\gamma\neq 0} \rightarrow 1$, with an approximate level $\alpha$ given by $\beta_n(0) \rightarrow  \alpha$.

%===============================================================================
\section{Comparative study via simulation}
\label{SecSims}

After dealing with the consistency and asymptotic distribution of the  general endpoint estimator $\hat{x}^F$ defined in \eqref{EPEst}, we are now ready to find out how these properties carry over to the finite sample setting. The  finite sample properties  of the test statistic  $G_{n,k}^*(0)$, defined in \eqref{Gnkg}, are also investigated,  assessing how it performs at either discerning the presence of a heavy-tailed model (with d.f. $F\in{\mathcal D}(G_\gamma)_{\gamma>0}$), or at detecting a short-tailed model (with $F\in{\mathcal D}(G_\gamma)_{\gamma<0}$). Consistency of the estimator (cf. Proposition \ref{PropCons}) justifies our belief that a larger sample leads to more accurate estimation  about the true value $x^F$, whereas consistency of the test (cf. Theorem \ref{Thm.Test}) connects a larger sample with a more powerful testing procedure. Of course how large ``sufficiently large" is, in terms of both $k$ and $n$, depends on the particular circumstances. The number of upper order statistics $k^*$ (yet to be determined) can be viewed as the effective sample size for extrapolation beyond the range of the available observations. In particular, if $k^*$ is too small, then the endpoint estimator tends to have a large variance, whereas if $k^*$ is too large, then the bias tends to dominate. This typical feature will be reflected in our simulation results. We argue comparison with other well-known endpoint estimators by means of their estimated absolute bias and mean squared errors.

To this end, we have generated $N=300$ samples from each of the four models listed below, taken as key examples:

\begin{itemize}
  \item \textbf{Model 1}, with d.f. $F_1(x)=1-\left[1+(-x)^{-\tau_1}\right]^{-\tau_2}$, $x<0$, $\tau_1,\tau_2>0$.   The EVI is $\gamma= -1/(\tau_1 \tau_2)$ and the endpoint $x^{F_1}=0$.\\
  \item \textbf{Model 2}, with d.f. $F_2(x)=1-\int_{-\infty}^{\log(1-1/x)} \lambda^2 t e^{-\lambda t}  dt$, $x<0$, $\lambda>0$.   The EVI is $\gamma=-1/\lambda$ and the endpoint $x^{F_2}=0$. Moreover, $X\id-1/(e^Z-1)$, where $Z$ is Gamma$(shape=2,rate=\lambda)$ distributed.\\
  \item \textbf{Model 3}, with d.f. $F_3(x)=1-\left[1+(\frac{1}{x}-1)^{-\tau_1}\right]^{-\tau_2}$, $x\in (0,1)$, $\tau_1,\tau_2>0$.   The EVI is $\gamma=-1/(\tau_1 \tau_2)$ and the endpoint is $x^{F_3}=1$.\\
\item \textbf{Model 4}, with d.f. $F_4(x)=1-(1-x)^{-1/\gamma}$,  $x\in (0,1)$, $\gamma<0$.   The EVI is $\gamma$ and the endpoint is $x^{F_4}=1$. This corresponds to a $Beta(1,-1/\gamma)$ model.
\end{itemize}
Each one of these models satisfies the main assumption that $F\in \mathcal{D}(G_{\gamma})$, for some EVI $\gamma<0$, which immediately entails a finite right endpoint $x^F$.

Models 1, 2 and 3 are the same ones as in \cite{GGS:11, GGS:12}, although these works only tackle the EVI equal to $-1$. Model 4 is a  Beta distribution parameterized in $\gamma<0$.  At the present stage we are interested in studying the exact performance of the general endpoint estimator  $\hat{x}^F$ for different ranges of the negative EVI. The extreme value index $\gamma$ is therefore a design parameter in the present simulation study, albeit under the restriction to $\gamma<0$. The second design parameter is of course the true right endpoint $x^F$, thus assumed finite. A number of combinations between model and design factors are assigned in order to obtain distinct values of the negative EVI, particularly
$\gamma= -1/2, -1/5$, together with two possibilities for the right endpoint, $x^F=0$ and $x^F=1$. The case $\gamma =0 $ has been extensively studied in \citet{FAN:14}, thus being obviated in the present setting.

\subsection{Endpoint estimation}

The finite sample performance of the general  estimator (notation: FAN) is here compared  with the na\"ive maximum estimator  $X_{n,n}$ (notation: MAX) and with the estimator $\hat{x}^*$ (notation: MOM.inv) introduced in equation (2.21) from \citet{FHP:03}:
\begin{equation}\label{mom_est}
\hat x^*:=X_{n-k,n}-\frac{\hat{a}(n/k)}{\hat{\gamma}_{n,k}^{-} }.
 \end{equation}
We note that the above estimator evolves from \eqref{EPClass} by using the consistent estimators for the EVI and scale function, respectively defined as
 \begin{equation}\label{EVI-inv}
\hat{\gamma}_{n,k}^{-} := 1-\frac{1}{2} \biggl\{1-\frac{\left(N_{n,k}^{(1)}\right)^2}{N_{n,k}^{(2)}}\biggr\}^{-1}
 \end{equation}
and
 \begin{equation}\label{ank-inv}
\hat{a}(n/k):= N_{n,k}^{(1)}(1-\hat{\gamma}_{n,k}^{-} ),
 \end{equation}
where
  \begin{equation}\label{moms-inv}
 N_{n,k}^{(r)}:=\frac{1}{k}\sum_{i=0}^{k-1}  \left( X_{n-i,n}-X_{n-k,n}\right) ^r, \quad r=1,2.
 \end{equation}
Note that EVI estimator \eqref{EVI-inv} is shift and scale invariant.
The class \eqref{EPClass} of estimators for the finite right endpoint has already been applied to a lifespan study by \citet{Aarssen:94} and \citet{FHP:03}, under the assumptions of a finite endpoint and  that the EVI lies between $-1/2$ and $0$. Another method for estimating the right endpoint for negative EVI is
via modeling the exceedances over a certain high threshold exactly by a  Generalized Pareto distribution (GPD). The result underpinning this parametric approach establishes that $F \in \mathcal{D}(G_{\gamma}) $  is equivalent to the relation
\begin{equation}\label{GP_POT}
\lim_{t\uparrow x^F} \sup_{0<x<x^F-t}\left| \frac{F(t+x)-F(t)}{1-F(t)} - H_\gamma \Bigl(\frac{x}{\sigma_t}\Bigr)\right|=0,
\end{equation}
where $H_\gamma(x):= 1+ \log G_\gamma(x)$ is the GPD \citep[see e.g.][]{deHaanFerreira:06}. In particular, if $\gamma=0$, the GPD reduces to the exponential  d.f. $H_0(x):= 1-\exp(-x),\,x\geq 0$.
In brief, condition \eqref{GP_POT} states that $F\in \mc{D}(G_\gamma)$ if and only if the excesses $Y:= X - t$ above a high threshold $t$ are asymptotically Generalized Pareto distributed.
Next to introducing  the class of GP distributions, relation \eqref{GP_POT} also enables to step away from the max-domain of attraction, towards the actual fit of the GPD to the sample excesses, providing a natural entry point to the POT approach. Once selected a \emph{high threshold} $t$, the POT method deems the shape parameter $\gamma \in \real$ (analogue to the EVI) and scale parameter $\sigma_t>0$ (which ultimately accommodates the influence of the threshold $t$) as the two indices characterizing the excess distribution function over $t$. Then we can proceed via maximum likelihood  (ML), the methodology at the core of the POTML.GPD procedures.  We note that, if $\gamma<0$, the parametric GPD fit corresponds to modeling the  exceedances $X$ over $t$  by a Beta distribution with finite right endpoint estimated by $\hat x^F_{POT}= t-\hat \sigma^{ML}/\hat \gamma^{ML}$, with $\hat \gamma^{ML}<0, \, \hat\sigma^{ML}>0$. Of course, in the case $\gamma=0$, a finite right endpoint is not allowed while fitting the exponential distribution.
References about the POTML.GPD approach are the seminal works by \citet{smith:87} and \cite{DavSmith:90}.

In the  semi-parametric setting, i.e. while working in the domain of attraction rather than dealing with the limiting distribution itself, the  upper intermediate
o.s. $X_{n-k,n}$ plays the role  of the \emph{high threshold} $t$. For the asymptotic properties of the POTML estimator of the shape parameter $\gamma$ under a semi-parametric approach, see e.g. \citet{drees:2004}, \citet{LiPeng:10} and \citet{Zhou:10}.

In case $F\in \mc{D}(G_\gamma)$ with $\gamma <0$, we have mentioned before endpoint estimators arising from the class \eqref{EPClass}.
Now, let
$\{ Y_i:=X_{n-i+1,n} - X_{n-k,n}\}_{i=1}^k $ be the excesses above the supposedly high random threshold $X_{n-k,n}$.  Furthermore, and building on the relation above, the sample  excesses $\{Y_i\}_{i=1}^n$ are assumed to follow a GPD. Then, the ML estimator $(\hat\sigma^{ML}_{k},\hat\gamma^{ML}_{k})$ can be worked out as the solution of
\begin{equation*}
\argmax_{\gamma<0, \,\sigma>0}
 \prod_{i=1}^k h_\gamma(Y_i/\sigma)/\sigma=\argmax_{\gamma<0, \,\sigma>0}  \prod_{i=1}^k \left(1+\frac{\gamma}{\sigma} Y_i\right)^{-(1/\gamma+1)} \sigma^{-k},
\end{equation*}
with $h_\gamma(x):=\frac{\partial}{\partial x} H_\gamma(x)$ \citep[see p.19 of][for a detailed explanation]{deHaanFerreira:06}. The POTML.GPD estimator of the right endpoint is defined as
\begin{equation}\label{POTML}
 \hat x^F_{ML} := X_{n-k,n}-\hat\sigma^{ML}_{k}/\hat\gamma^{ML}_{k}\,,\,\mbox{ with }\, \hat\gamma^{ML}_{k}<0 \,\mbox{ and } \, \hat\sigma^{ML}_{k}>0,
\end{equation}
showing a close similarity with  the semi-parametric class \eqref{EPClass} of right endpoint estimators. We refer to section 4.5.1 of \citet{deHaanFerreira:06}  and \citet{QiPeng:09} for the semi-parametric handling of \eqref{POTML}.

The POTML.GPD endpoint estimator relies on the shift and scale invariant ML-estimator of the shape parameter $\gamma<0$, a restriction strictly to ensure that \eqref{POTML} is a valid endpoint estimator. However there is no explicit formula for the ML-estimator. An accumulating literature has pointed out this disadvantage. Maximization of the log-likelihood, reparameterized in $(\tau:=\gamma/\sigma, \gamma)$, has been discussed in \citet{Grimshaw:93}. Although theoretically well determined, even when $\gamma \uparrow 0$, the non-convergence to a ML-solution can be an issue when $\gamma$ is close to zero. There are also irregular cases which may compromise the practical applicability of ML. Theoretical and numerical accounts of these issues can be found in \citet{CastilloDaoudi:09} and \citet{CastilloSerra:15} and references therein.

Inspired by the numeric examples in \cite{GGS:11, GGS:12}, we have generated $N=300$ replicates of a random sample with size $n = 1000$ and computed the average $L^1$-error given by
 \begin{equation*}
   E(k^*):=\frac{1}{N}\sum_{j=1}^{N}|\varepsilon(j,k^*)|, \quad \mbox{where} \quad \varepsilon(j,k^*):={ \hat x_{k^*}}(j) -x^F, \,\,\,k^*\leq n,
 \end{equation*}
where $x_{k^*}(j)$ denotes the endpoint estimator evaluated at the $j$-th replicate, for every $k^*$.

We borrow models 1-3 from \cite{GGS:11, GGS:12} and their performance measures but we will not proceed with their proposals for endpoint estimation. Unlike their high order moment estimators, none of the endpoint estimators adopted in the present simulation study (MAX, FAN, MOM.inv, and POTML.GPD) require the knowledge of the original sample size $n$, since these rely on a certain number $k^*$  of top  o.s. only. For the purpose of simplicity, the number $k^*$ will be viewed as the effective sample size.
In the sequel, the na\"ive estimator MAX is attached to $k^*=1$ since it coincides with the first top o.s.; the POTML.GPD and MOM.inv  endpoint estimators are functions of the $k^*=k+1$ upper o.s.; finally, the FAN estimator requires  $k^*=2k$ top observations. We can also compute the ``optimal" values of $k^*$ in the sense of minimizing the average $L^1$-error, i.e. $k_0^*:=\arg\min\{ E(k^*),\, k^*\leq n\}$.
 Since the MAX entails $k^*=1$, the associated function $E$ is constant
 and this optimality criterium has no effect on the na\"ive estimator.  Table \ref{table1} displays the simulation results, where we have considered parameter combinations with respect to $\gamma=-0.5, -0.2$ and $x^F=0,\,1$. The POTML.GPD endpoint estimates were found by maximizing the log-likelihood over $\gamma <0$.
\begin{table}
  \caption{\footnotesize Average $L^1$-errors. The lowest values appear in bold.}\label{table1}
  \vspace{0.5cm}\centering
   \begin{tabular}{llcccc}
    \hline
    \\
    Model &   \hspace{0.5cm} & MAX & FAN   & MOM.inv & POTML.GP\\
    \\
    \hline
    \\
    Model 1  $(x^F=0)$ \\
    $(\tau_1,\tau_2)=(2,1)$ &  \hspace{0.3cm}& 0.028 & \textbf{0.013} & 0.024 &0.019 \\
    $(\tau_1,\tau_2)=(5,1)$ &  \hspace{0.3cm}& 0.231 & \textbf{0.041} & 0.145 &0.128 \\
    \\
    Model 2 $(x^F=0)$\\
    $\lambda=2$ &  \hspace{0.3cm}&  0.009 & \textbf{0.005} & 0.009 &0.008\\
    $\lambda=5$ &  \hspace{0.3cm}&  0.173 & \textbf{0.044} & 0.133 &0.133\\
    \\
    Model 3 $(x^F=1)$\\
    $(\tau_1,\tau_2)=(2,1)$ &  \hspace{0.3cm}& 0.027 & \textbf{0.012} & 0.022 & 0.014 \\
    $(\tau_1,\tau_2)=(5,1)$ &  \hspace{0.3cm}& 0.187 & 0.129 & 0.120 & \textbf{0.095} \\
    \\
    Model 4 $(x^F=1)$\\
    $\gamma=-1/2$ &  \hspace{0.3cm}&  0.029  & \textbf{0.014 }& 0.029 &\textbf{0.014}\\
    $\gamma=-1/5$ &  \hspace{0.3cm}&  0.234 & 0.171 &  0.103 & \textbf{0.085}\\
    \hline
  \end{tabular}
\end{table}
\noindent

\begin{figure}
  \centering
  \includegraphics[width=0.4\textwidth]{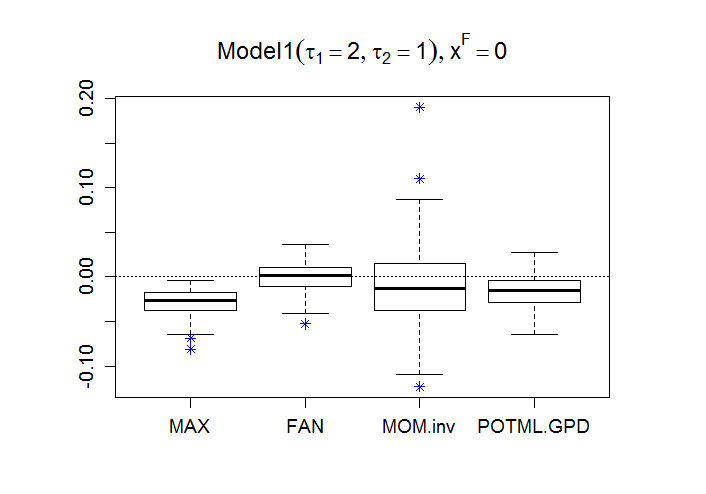}
  \includegraphics[width=0.40\textwidth]{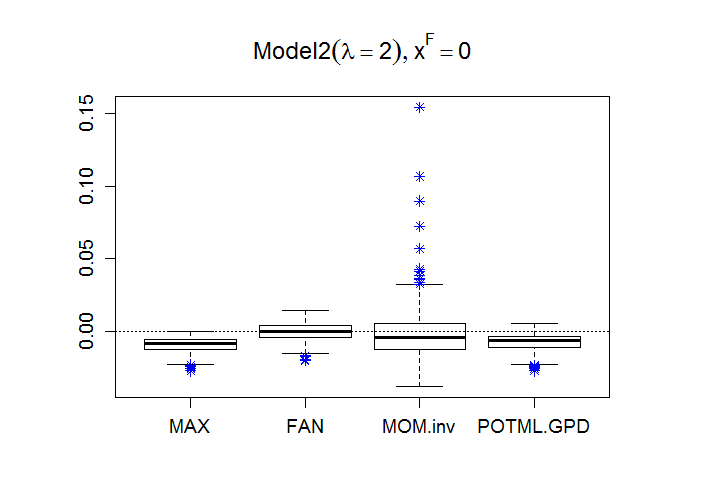} \\
       \includegraphics[width=0.4\textwidth]{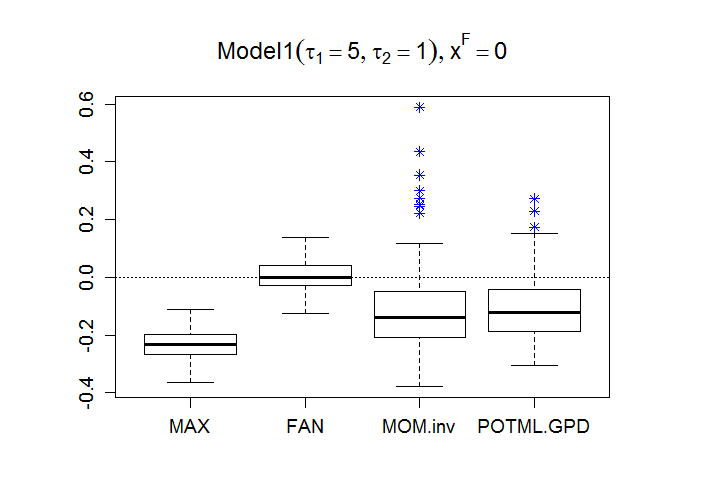}
       \includegraphics[width=0.4\textwidth]{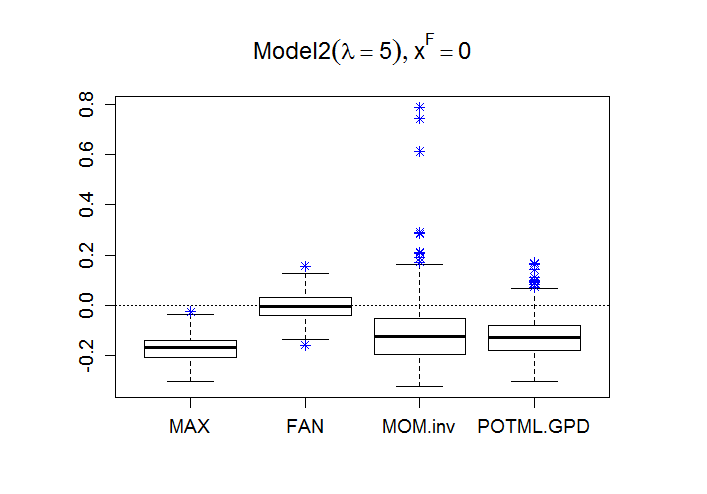}
       \caption{\footnotesize  Boxplots of the optimal bias $\varepsilon(j,k^*_0),\,\,j=1,\ldots, N$, with $N=300$. Endpoint estimates are drawn from Model 1 (\emph{left}) and Model 2 (\emph{right}), with true value $x^F=x^{F_1}=x^{F_2}=0$.}
      \label{boxes12}
\end{figure}

\begin{figure}
  \centering
  \includegraphics[width=0.4\textwidth]{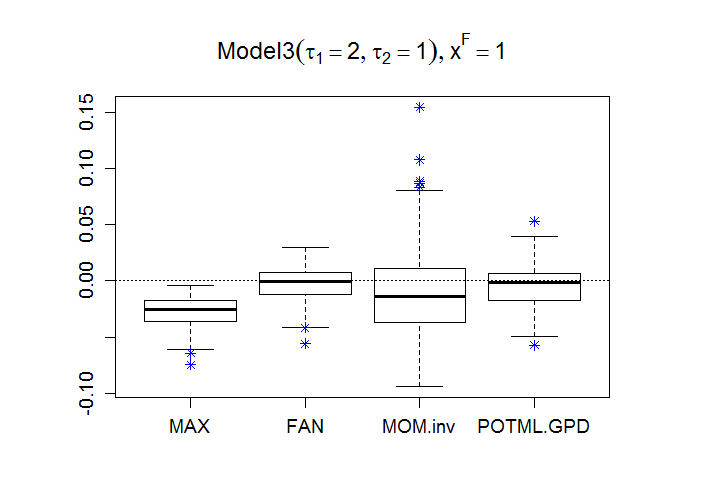}
   \includegraphics[width=0.4\textwidth]{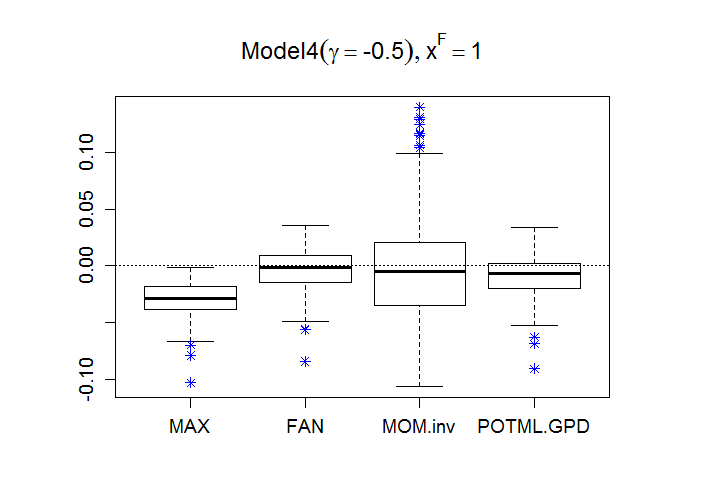}\\
       \includegraphics[width=0.4\textwidth]{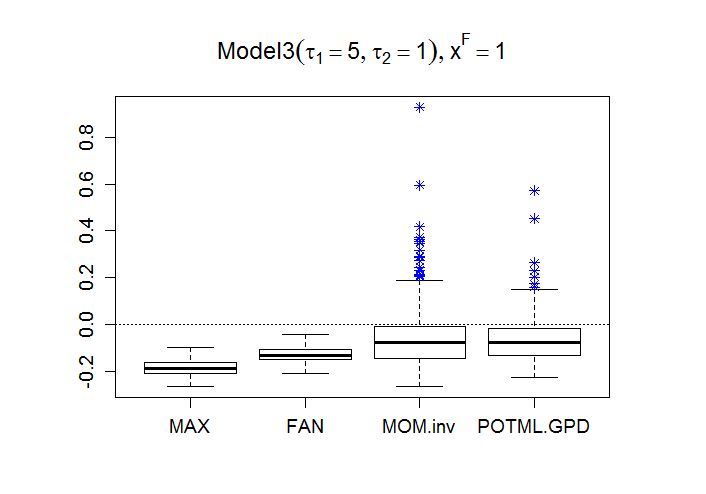}
        \includegraphics[width=0.4\textwidth]{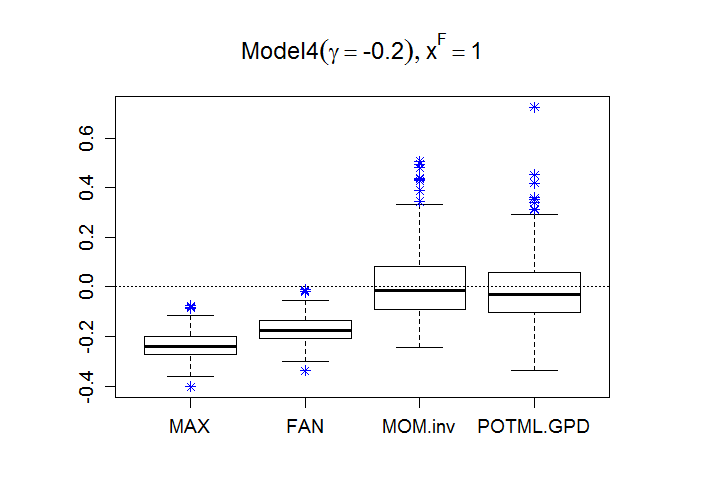}
        \caption{\footnotesize  Boxplots of the optimal bias $\varepsilon(j,k^*_0),\,\,j=1,\ldots, N$, with $N=300$. Endpoint estimates are drawn from Model 3 (\emph{left}) and Model 4 (\emph{right}), with true value $x^F=x^{F_3}=x^{F_4}=1$.}
      \label{boxes34}
\end{figure}

The relative performance of the adopted endpoint estimators on the ``optimal"  $k_0^*$ is depicted in the  box-plots of Figures \ref{boxes12} and \ref{boxes34}, in terms of their associated errors $\varepsilon(j,k^*_0),\,\,j=1,\ldots, N$, with $N=300$. Apart from the obvious conclusion that the MAX tends to underestimate the true value of the endpoint $x^F$,  we find that the POTML.GPD, MOM.inv and FAN estimators have distinct behaviors with respect to the optimal levels $k^*_0$. In particular, FAN estimates are not so spread out as the ones returned by POTML.GPD or by the MOM.inv endpoint estimators, the latter showing larger variability.

\begin{figure}
  \centering
 \includegraphics[width=0.4\textwidth]{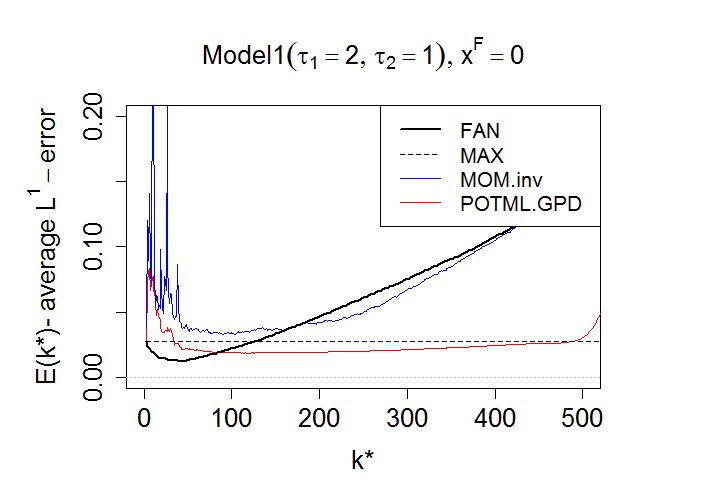}
   \includegraphics[width=0.4\textwidth]{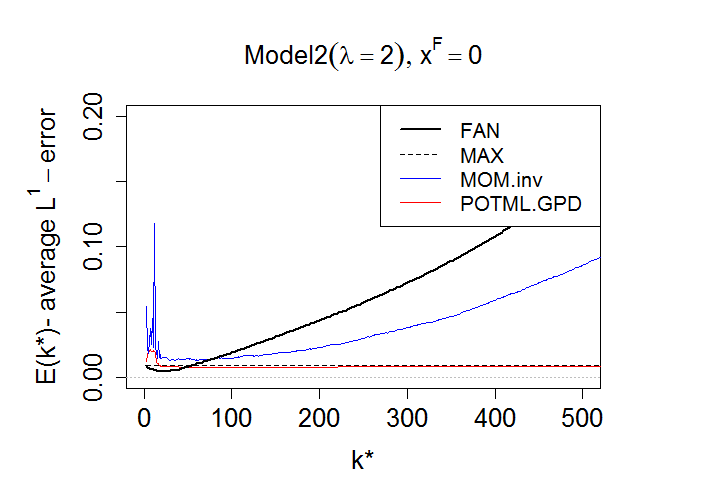}\\
  \includegraphics[width=0.4\textwidth]{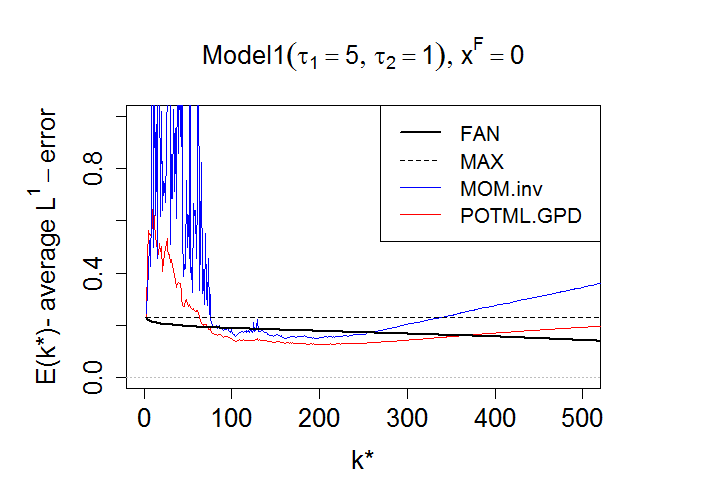}
  \includegraphics[width=0.4\textwidth]{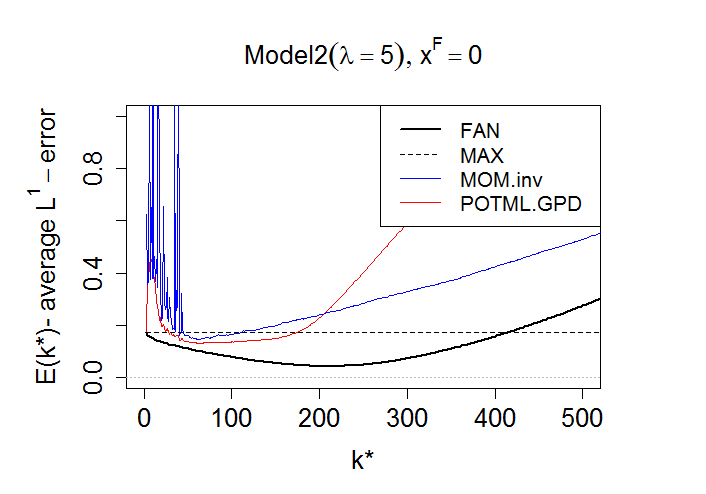}
  \caption{\footnotesize Average $L^1$-error, $E(k^*)$, plotted against $k^*\leq n/2$.  Endpoint estimates are drawn from Model 1 (\emph{left}) and Model 2 (\emph{right}), with true value $x^F=x^{F_1}=x^{F_2}=0$.}
      \label{E12}
\end{figure}

\begin{figure}
  \centering
  \includegraphics[width=0.4\textwidth]{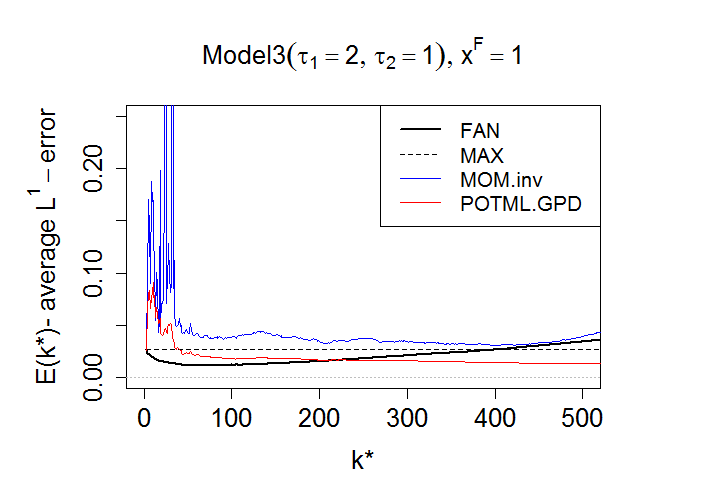}
 \includegraphics[width=0.4\textwidth]{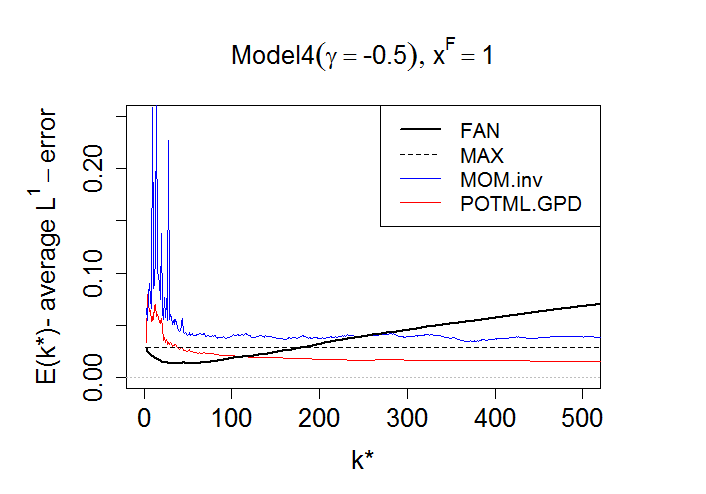}\\
  \includegraphics[width=0.4\textwidth]{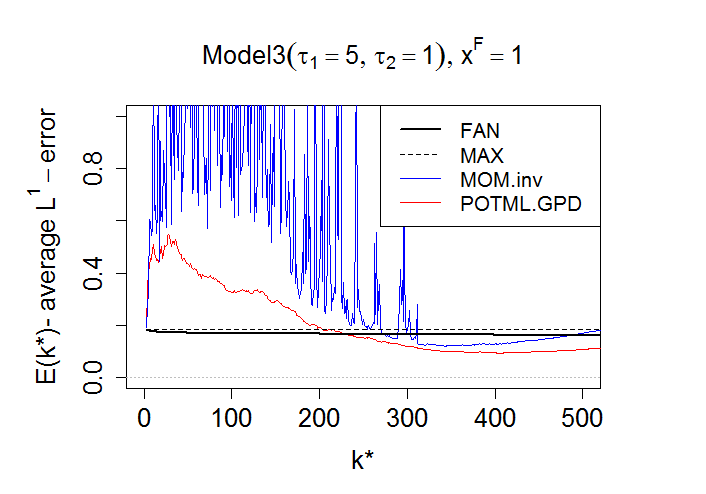}
  \includegraphics[width=0.4\textwidth]{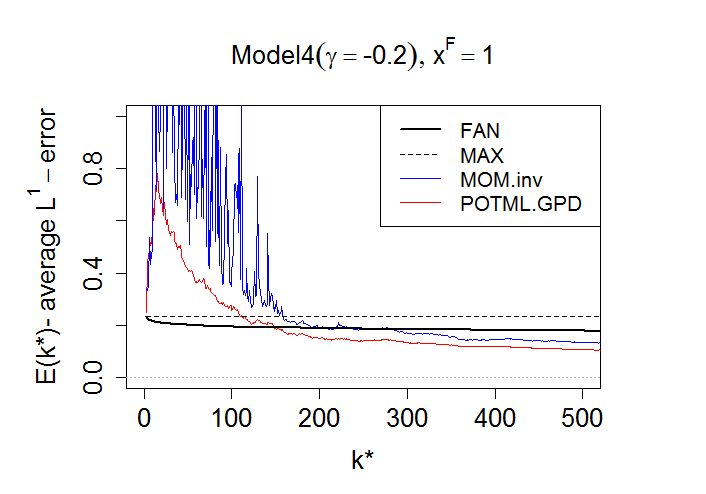}
    \caption{\footnotesize Average $L^1$-error, $E(k^*)$, plotted against $k^*\leq n/2$. Endpoint estimates are drawn from Model 3 (\emph{left}) and Model 4 (\emph{right}), with true value $x^F=x^{F_3}=x^{F_4}=1$.}
       \label{E34}
\end{figure}

\begin{figure}
  \centering
   \includegraphics[width=0.4\textwidth]{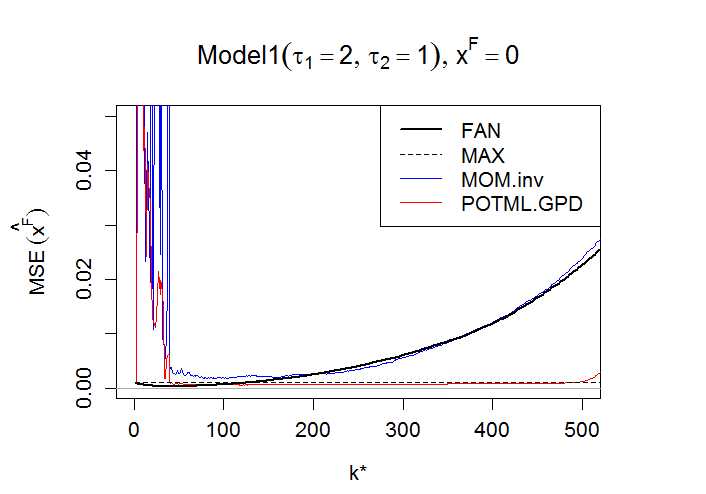}
  \includegraphics[width=0.4\textwidth]{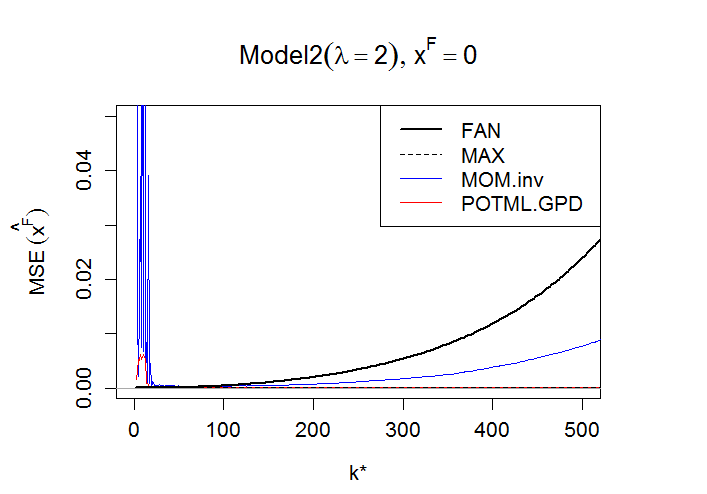}\\
 \includegraphics[width=0.4\textwidth]{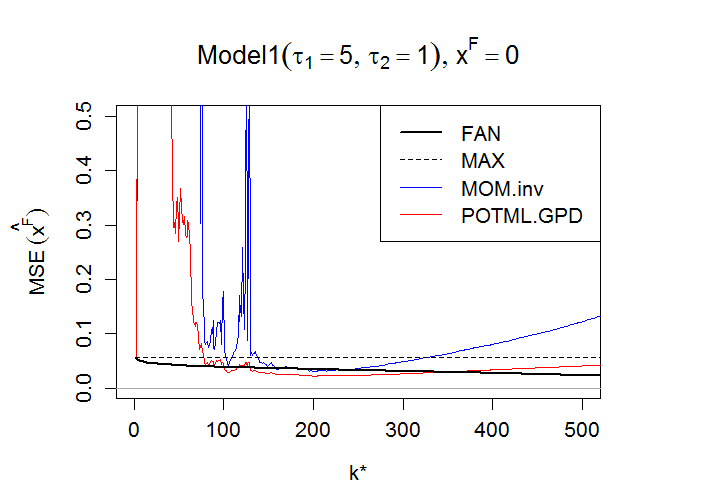}
  \includegraphics[width=0.4\textwidth]{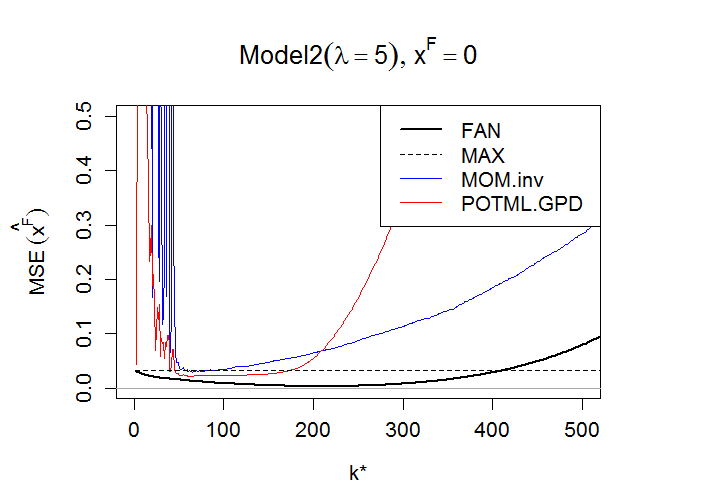}
    \caption{\footnotesize Mean squared error (MSE) as function of the number of upper o.s. $k^*$,  $k^*\leq n/2$. Endpoint estimates are drawn from Model 1 (\emph{left}) and Model 2 (\emph{right}), with true value $x^F=x^{F_1}=x^{F_2}=0$.}
      \label{mse12}
\end{figure}

\begin{figure}
    \centering
   \includegraphics[width=0.4\textwidth]{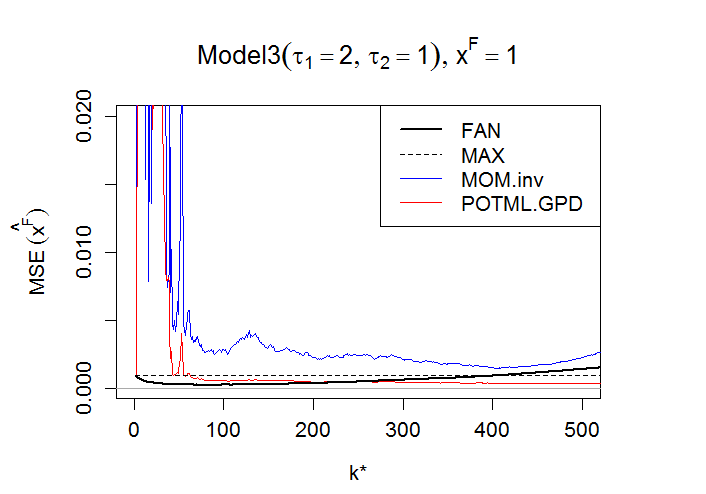}
\includegraphics[width=0.4\textwidth]{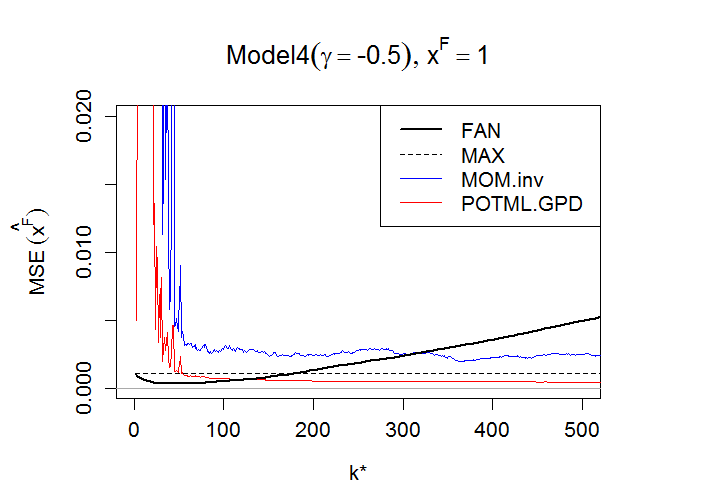}\\
     \includegraphics[width=0.4\textwidth]{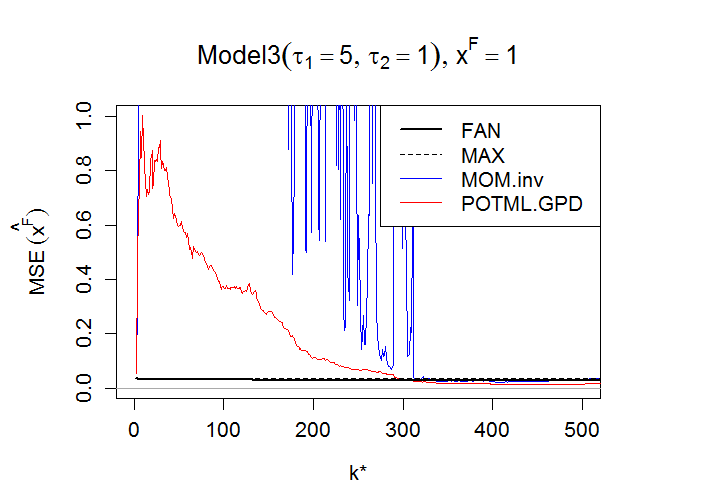}
  \includegraphics[width=0.4\textwidth]{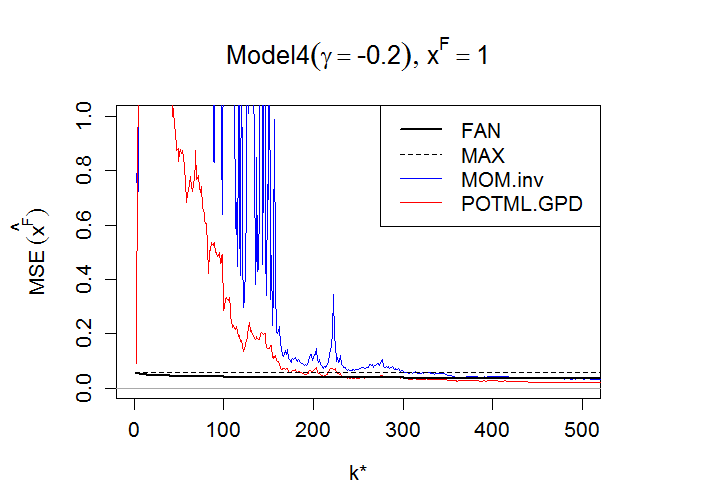}
    \caption{\footnotesize Mean squared error (MSE) as function of the number of upper  o.s. $k^*$,   $k^*\leq n/2$. Endpoint estimates are drawn from Model 3 (\emph{left}) and Model 4 (\emph{right}), with true value $x^F=x^{F_3}=x^{F_4}=1$.}
      \label{mse34}
\end{figure}

Figures \ref{E12} and \ref{E34} display the plain average $L^1$-error (i.e. without the optimality assessment) against  the number $k^*$ of upper o.s. used in the corresponding estimation process. The pertaining mean squared errors (MSE) are depicted in Figures \ref{mse12} and \ref{mse34}, respectively. The four  models here addressed are set with EVI$=-1/2,-1/5$. From these results, it is clear that the MOM.inv and POTML.GPD endpoint estimators are very unstable for small values of $k^*$ ($k^*\leq 100$), contrasting with the small variance of the proposed FAN estimator along the entire trajectory.  On the other hand, both MOM.inv and FAN estimators  show increasing $L^1$-error with increasing of $k^*$, a common feature to extreme semi-parametric estimators. The FAN estimator seems to perform best in those regions where other estimators exhibit high volatility, which may range from small to moderately large values of $k^*$. This feature is more severe when $\gamma=-0.2$ (see second row of Figures \ref{E12} and \ref{E34}), where the instability persists until an impressive $k^*=300$ is reached. Once attained a plateau of stability, the POTML.GPD tends to perform very well in general. The best way to apply MOM.inv seems to dwell in a precise choice of $k^*$, which should be selected at the very end of the very erratic path, just before bias sets in.
The general endpoint estimator (FAN) tends to return values with a low average $L^1$-error and low MSE. In fact, Figures \ref{E12}  up to \ref{mse34} not only provide us with a snapshot for this specific choice of EVI values ($-1/2$ and $-1/5$), but also allow to foresee the estimates behavior with respect to other EVIs  in between, once we screen the plots from the top to the bottom in each Figure. The boxplots in Figures \ref{boxes12} and \ref{boxes34} already suggested this possibility: the outliers marked in these boxplots seem to move from lower to larger values of optimal bias $\varepsilon(\cdot,k^*_0)$ as we progress on increasing EVI.
%==============================================================
 
%============ Boxplots========================

 %%%%%%%%%%%%%%%%%%%%%%%% E()%%%%%%%%%%%%%%%%%%%%%%%%%%
%%%%%%%%%%%%%%%%%%%%%%%% MSE()%%%%%%%%%%%%%%%%%%%%%%%%%%

Altogether, the general endpoint estimator FAN seems to be an improvement to the na\"ive MAX estimator and tends to surpass the MOM.inv and POTML.GPD estimators, by delivering low biased estimates quite often, while showing a low variance component.
This is particularly true for a EVI close to zero, as it would be expected from the one estimator primarily tailored to tackle endpoint estimation in the Gumbel max-domain of attraction \citep[cf.][]{FAN:14}.
Furthermore, the FAN estimator seems to work remarkably well under a fairly negative EVI, considering that this is a general estimator which does not accommodate any specific information about the true value of the EVI. The overall performance of the MOM.inv and  POTML.GPD  endpoint estimators is clearly damaged by their  large variance in the top of the sample. It is worthy to notice that the presented estimation procedure $\hat x^F$ acts as a complement to numerical POTML methods for endpoint estimation, the latter only available for strict negative shape parameter $\gamma$; in contrast, $\hat x^F$ presents an explicit simple expression, unifying the estimation method to distributions with non-positive EVI, $\gamma\leq 0$.

%----------------------------------------
\subsection{A testing procedure built on the general endpoint estimator} \label{newTest}

This section concerns the finite sample performance of  $G_{n,k}$, presented in Theorem \ref{Thm.Test}, as a convenient tool  for either discarding  heavy-tailed models or  for detecting short-tailed models $F\in{\mathcal D}(G_\gamma)$, $\gamma<0$.
 One grounding result in this respect is that all the semi-parametric endpoint estimators we are adopting, are consistent under the assumption that
 $k=k_n$ is an intermediate sequence of positive integers, i.e. $k=k_n \rightarrow \infty$ and $k_n/n \rightarrow 0$, as $n\rightarrow \infty$. The testing procedures we wish to apply also bear on this usual assumption in statistics of extremes.
There are many proposals for testing procedures aiming at the selection of a suitable max-domain of attraction. For a wide view on this topic, we refer the surveys on testing about extreme values conditions available in \citet{HuslerPeng:08} and \citet{ NevesFA:08}. We recall that  EVI estimation  is not a requirement for the
general endpoint estimation defined \eqref{EPEst} and emphasize that both Weibull and Gumbel domain are allowed. Thus, for the time being, we will rely on testing procedures which do not require external estimation of the EVI. We will compare the new test statistic $G_{n,k}$, defined in \eqref{G}, with  the Ratio and Greenwood statistics introduced in \citet{NevesPFA:06} and \citet{NevesFA:07} for the one-side alternatives
\begin{eqnarray*}
H_0: F\in{\mathcal D}(G_0)\qquad vs   & &H_1: F\in{\mathcal D}(G_\gamma)_{\gamma>0} \\
&[\mbox{or  } &H_1': F\in{\mathcal D}(G_\gamma)_{\gamma<0}] \,.
\end{eqnarray*}
The Ratio(R) and Greenwood(Gr) statistics are defined as
\begin{equation*}
 R:=\frac{X_{n,n}-X_{n-k,n}}{N_1} \qquad Gr:=\frac{N_2}{\left(N_1\right)^2} \,\,,\\
\end{equation*}
with $N_j=N^{(j)}_{n,k}=\frac{1}{k}\sum_{i=0}^{k-1}  \left( X_{n-i,n}-X_{n-k,n}\right) ^j,\,\,j=1,2$.
 Under $H_0$,  the standardized version of Ratio statistic, $R^*:=R-\log k$, is asymptotically Gumbel, whereas the suitably normalized Greenwood  statistic, $Gr^*:=\sqrt{k/4}(Gr-2)$, is asymptotically standard normal.
Formally, approximated  $\alpha$-significant tests against the alternative $H_1$ [resp. $H_1'$] render the rejection regions $  R^* \geq\omega_{1-\alpha}$ [resp. $R^* \leq\omega_{\alpha}$] and $ Gr^* \geq z_{1-\alpha}$ [resp. $Gr^* \leq z_{\alpha}$], where $\omega_{\epsilon}:=\Lambda^{\leftarrow}(\epsilon)$ and $z_{\epsilon}:=\Phi ^{\leftarrow} (\epsilon)$. Here, $\Phi$ denotes the d.f. of the standard normal. Corresponding $p$-values of the test are  $p=1-\Lambda(g^*)$ against the heavy-tailed alternative, and  $p=\Lambda(g^*)$ against the short-tailed alternative, for  the observed values $r^*$  and $gr^*$  of the test statistics $R^*$ and $Gr^*$, respectively. The  approximated $p$-values against heavy-tailed alternatives $H_1$ [resp. short-tailed alternatives $H_1'$] are given by $1-\Lambda(r^*)$ and $1-\Phi(gr^*)$ [resp. $\Lambda(r^*)$ and $\Phi(gr^*)$], for an observed value $g^*:=g_{n,k}^*(0)$ of the test statistic $G_{n,k}^*(0)$.

\begin{figure}
    \centering
   \includegraphics[width=0.4\textwidth]{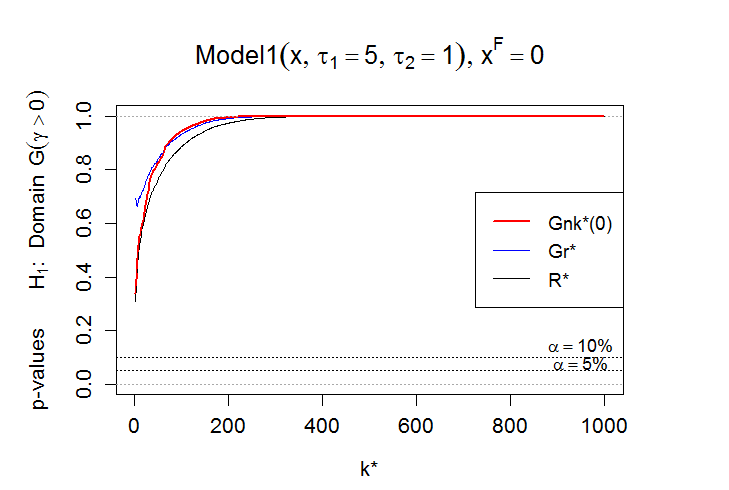}\hfil
     \includegraphics[width=0.4\textwidth]{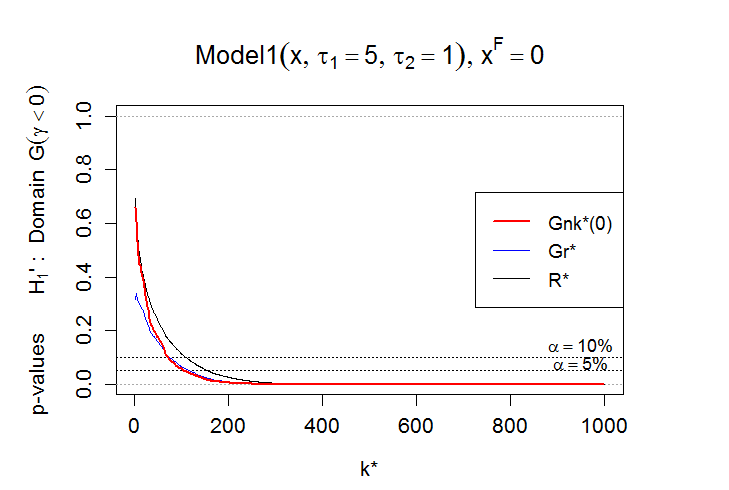}\\
        \includegraphics[width=0.4\textwidth]{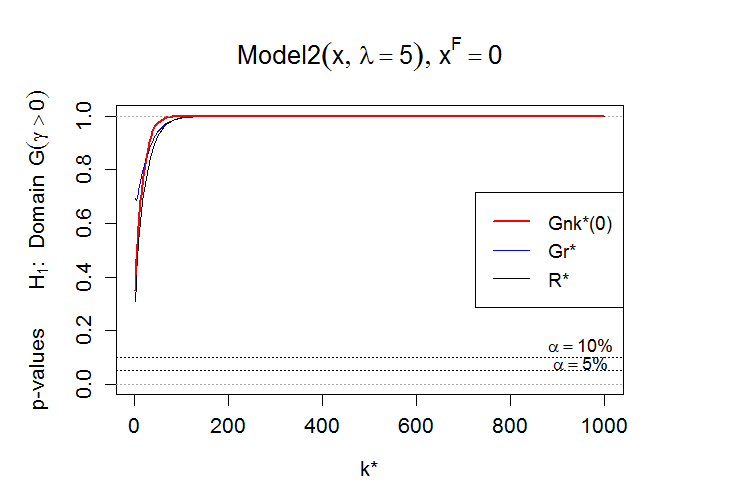}\hfil
     \includegraphics[width=0.4\textwidth]{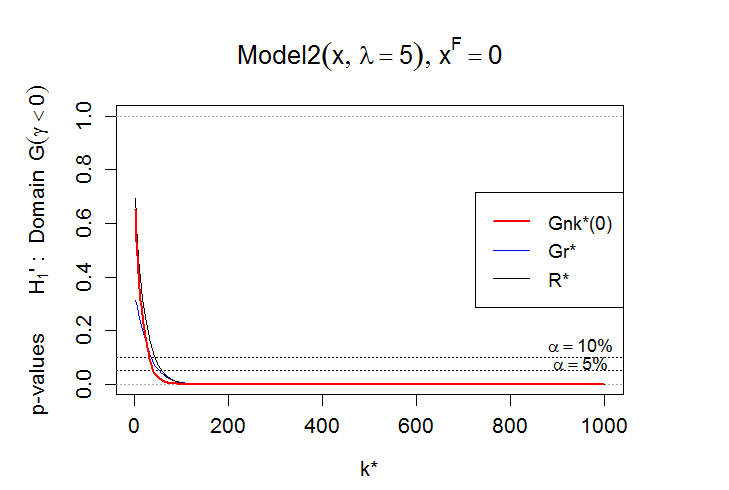}\\
    \caption{\footnotesize Average $p$-values of the simulated $G_{n,k}^*(0)$, $R^*$ and $Gr^*$ either against a heavy-tailed alternative (\emph{left}) or a short-tailed alternative (\emph{right}), with respect to $k^*\leq n \,=1000$.}
      \label{p-value1}
\end{figure}
\begin{figure}
    \centering
        \includegraphics[width=0.4\textwidth]{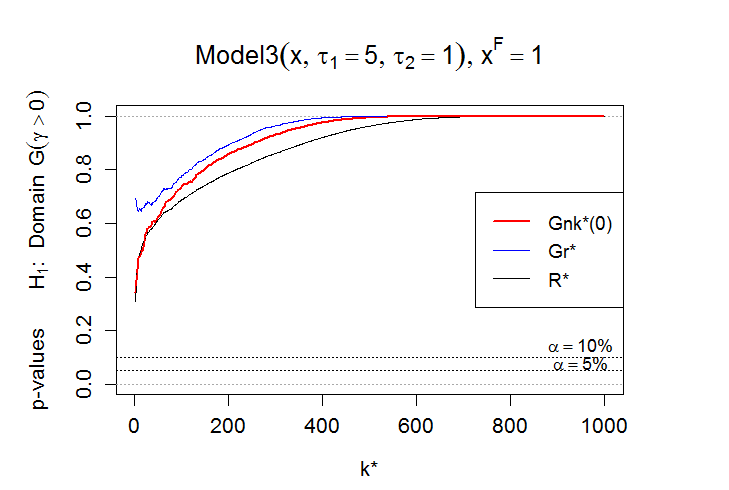}\hfil
     \includegraphics[width=0.4\textwidth]{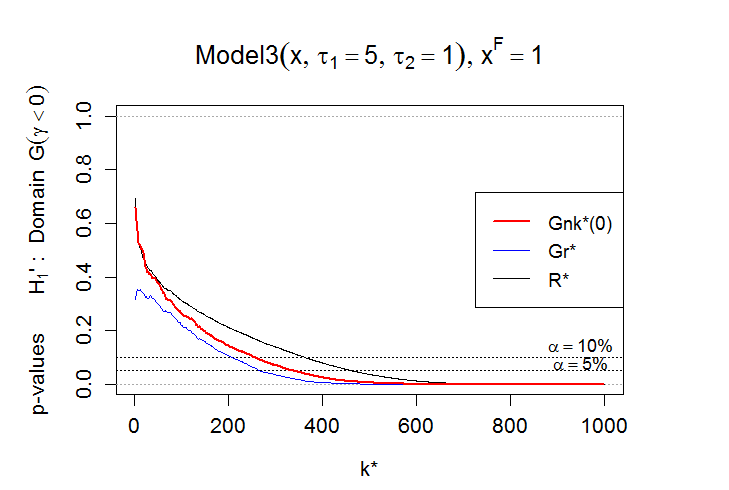}\\
       \includegraphics[width=0.4\textwidth]{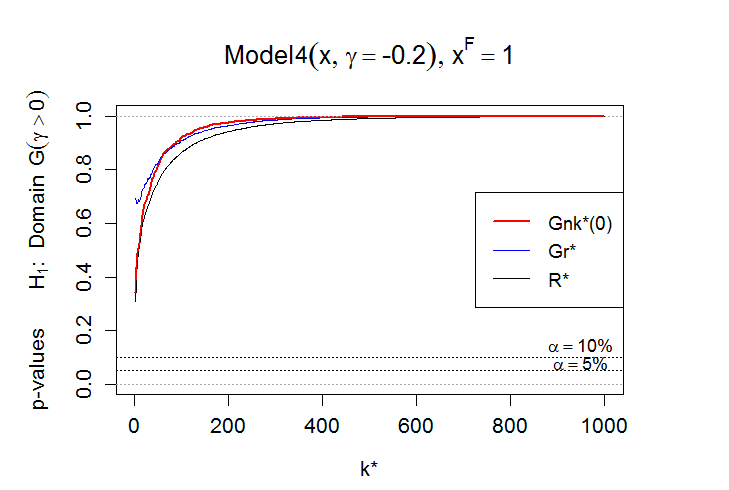}\hfil
     \includegraphics[width=0.4\textwidth]{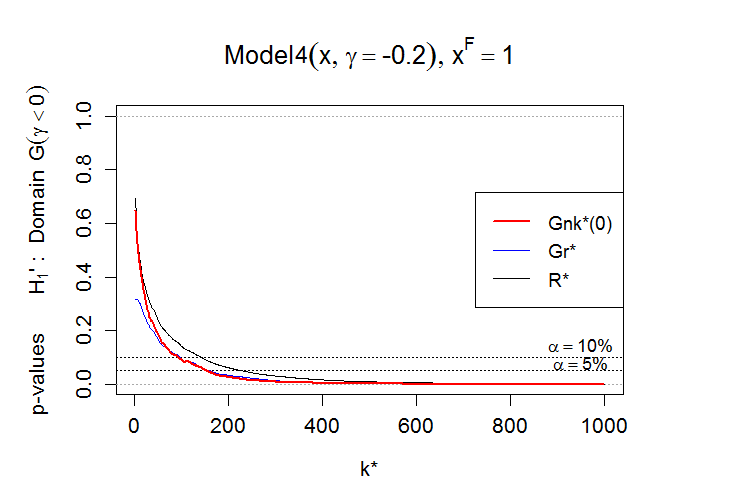}\\
      \caption{\footnotesize Average $p$-values of the simulated $G_{n,k}^*(0)$, $R^*$ and $Gr^*$ either against a heavy-tailed alternative (\emph{left}) or a short-tailed alternative (\emph{right}), with respect to $k^*\leq n \,=1000$.}
      \label{p-value2}
\end{figure}

 Figures \ref{p-value1}-\ref{p-value2} summarize the comparison between the performance of the test based on  $G_{n,k}^*(0)$ (cf. Theorem \ref{Thm.Test})  with the above mentioned tests $R^*$ and $Gr^*$. The simulations yield large  $p$-values  in connection with the heavy-tailed alternatives $H_1$, meaning that heavy-tailed distributions are likely to be detected by these tests. On the opposite side, the test are not so sharp against short-tailed alternatives in $H_1'$. The new test statistic  $G_{n,k}^*(0)$  rejects on smaller values of $k^*$ than the $R^*$ statistic, thus revealing more powerful than the Ratio-test.

The Greenwood test compares favourably to $G_{n,k}^*(0)$  in terms of power, against the heavy-tailed alternative.  However, it tends to be a more conservative test than the new proposal, often returning p-values much less than $5\%$.

%==================================================================

\section{Case study:  supercentenarian women lifespan}
\label{CaseStudy}

This section is devoted to the practical illustration of our methodology for statistical inference about the endpoint.
Our data set of oldest people comprises records of lifetimes in days of \emph{verified supercentenarians (women)}, with deaths in the time window 1986-2012. The  data set was extracted from  Table B of  Gerontology Research Group (GRG), as of January 1, 2014, merged with Tables C and E, as of June 29, 2015, available at \url{http://www.grg.org/Adams/Tables.htm}. Although the referred database includes lifespan records tracing back to 1903, these are often sparse and with a low average number of yearly records, which is not surprising since the GRG was only founded  in 1990. The later two years of records 2013-2014 are not yet closed. Therefore we settle with the $1272$ supercentenarian women lifetimes, recorded from 1986 to 2012, and corresponding to approximately $90\%$ of the total number of records since 1903.

The terms ``life expectancy" and ``lifespan" describe two entirely different concepts, although people tend to use these terms interchangeably. Life expectancy refers to the number of years a person is yet expected to live at any given age, based on the statistical average. Lifespan, on the other hand, refers to the maximum number of years that a person can potentially expect to live based on the greatest number of years anyone has lived. We are interested in the latter.

Formally, in gerontology literature, maximum lifespan potential (MLSP) is the operative definition for the verified age of the longest lived individual for a species \citep{OlshanskyCC:90} and, in this sense, can be viewed as a theoretical upper limit to lifetime.  The oldest documented age reached by any living individual is 122 years, meaning humans are said to have a MLSP of 122 years.
In Biology, theories of ageing are mainly divided into two groups: damage theories and program theories. According to damage theories, we age because our systems break down over time; so, if damage theories hold true, we can survive  longer by avoiding damaging our organism. Program theories consider that we age because there is an inbuilt mechanism that tells us to die; according to that, we cannot survive longer than the upper limit of longevity despite of our best efforts (see \citet{Hanayama:13}). \citet{Kaufmann:07} also discussed the issue of whether the right endpoint of the life span is infinite, for which they analyzed mortality data from West Germany. Their estimated shape parameter within the Generalized Pareto model is equal to $-0.08$ and the right endpoint of the estimated beta distribution is equal to 122 years, an estimate that fits well to the worldwide reported life span of the most famous record-holder, the Frenchwoman Jeanne-Louise Calment (122 years and 164 days) who was born in Arles on Feb. 21, 1875 and died in Arles on Aug. 4, 1997. However, \citet{Kaufmann:07} did not conclude categorically that human lifespan has a finite upper limit, arguing that by using the concept of penultimate distributions we can show that an infinite upper limit is well compatible with extreme value theory. They carry on pointing out a Beta distribution as a suitable model. \citet{Aarssen:94} analyzed lifespan data from the Netherlands using statistical methods under  the extreme value theory umbrella. \citet{Aarssen:94}  showed that there is a finite age limit, tackled with reasonable confidence bounds in the $113-124$ year span, a conclusion confined to the years of birth $1877-1881$ in the Netherlands.

Stephen Coles, a specialist in tracking human supercentenarians and co-founder of  the Supercentenarian Research Foundation (SRF), refers to the supercentenarians as ``the most extreme example of human longevity that we know about, the oldest old". In \citet{Coles:11}, the value 122 is referred to as the ``Calment Limit" for human longevity (what is designated here as the MLSP), which is supported on the fact that nobody has come even close that extreme age over the last 19 years. It is also mentioned that, from his research experience, ``supercentenarians had virtually nothing in common: they had different occupations, lifestyles, religions and so on, regardless the common factor of long-lived relatives." Also, according to \citet{Vaupel:11} ``the explosion in very long life has already begun", although by his perspective ``we cannot see much beyond 122."

Several authors have stated that despite of the increasing  ``life expectancy", the ``maximum human lifespan" has not much changed. According to \citet{TroenC:01} some biodemographic estimates predict that elimination of most of the major diseases  such as cancer, cardiovascular disease, and diabetes would add no more than 10 years to the average life expectancy, but would not affect MLSP \citep{Olshansky:90,TroenC:01}.  Other researchers go further enough to hypothesize that \emph{mortality will be compressed against that fixed upper limit} to life time  \citep[\emph{Compression Theory} by][]{Fries:80}. On the other hand,
\citet{WilmothR:03} argue  a possible  world trend in maximum lifespan, based on a long series of Swedish data.
Above all, there is still plenty of scope to assess significance of other covariates, like the negative impact of obesity and epidemic diseases on the rise in life expectancy  trends and the possible impact on the MLSP.

What researchers seem to agree on is the need for better data, since at present, there is insufficient data available on the extreme elderly population. We should keep in mind that age is often  misreported and at the time the centenarians (and \emph{supercentenarians}) were born, record keeping was less complete than it is nowadays.

With this  illustrative example of estimation of the ultimate lifespan by adopting the general endpoint estimator, it is not our aim to make conclusions for a specific cohort of individuals in time or space, nor any other type of serial studies. Instead, the interest will be  on the question of what are sensible bounds for the MLSP, at the current state of the art.\\

\begin{figure}
    \centering
      \includegraphics[width=0.7\textwidth]{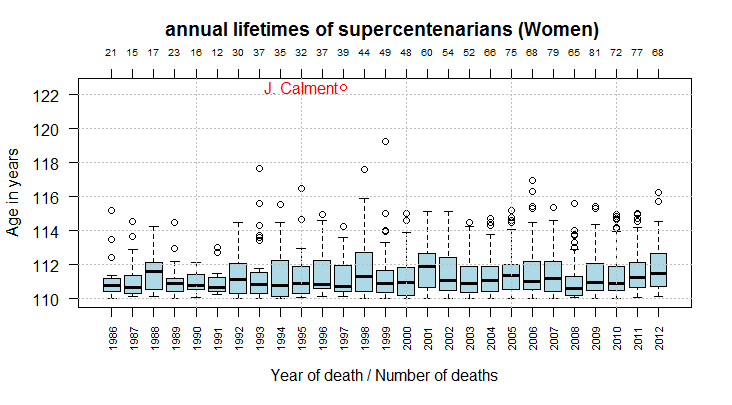}
     \caption{\footnotesize Comparative boxplots of the lifetimes of supercentenarians (Women) reported from 1986 to 2012.}
      \label{Box_Hist}
\end{figure}

\begin{figure}
    \centering
        \includegraphics[width=0.6\textwidth]{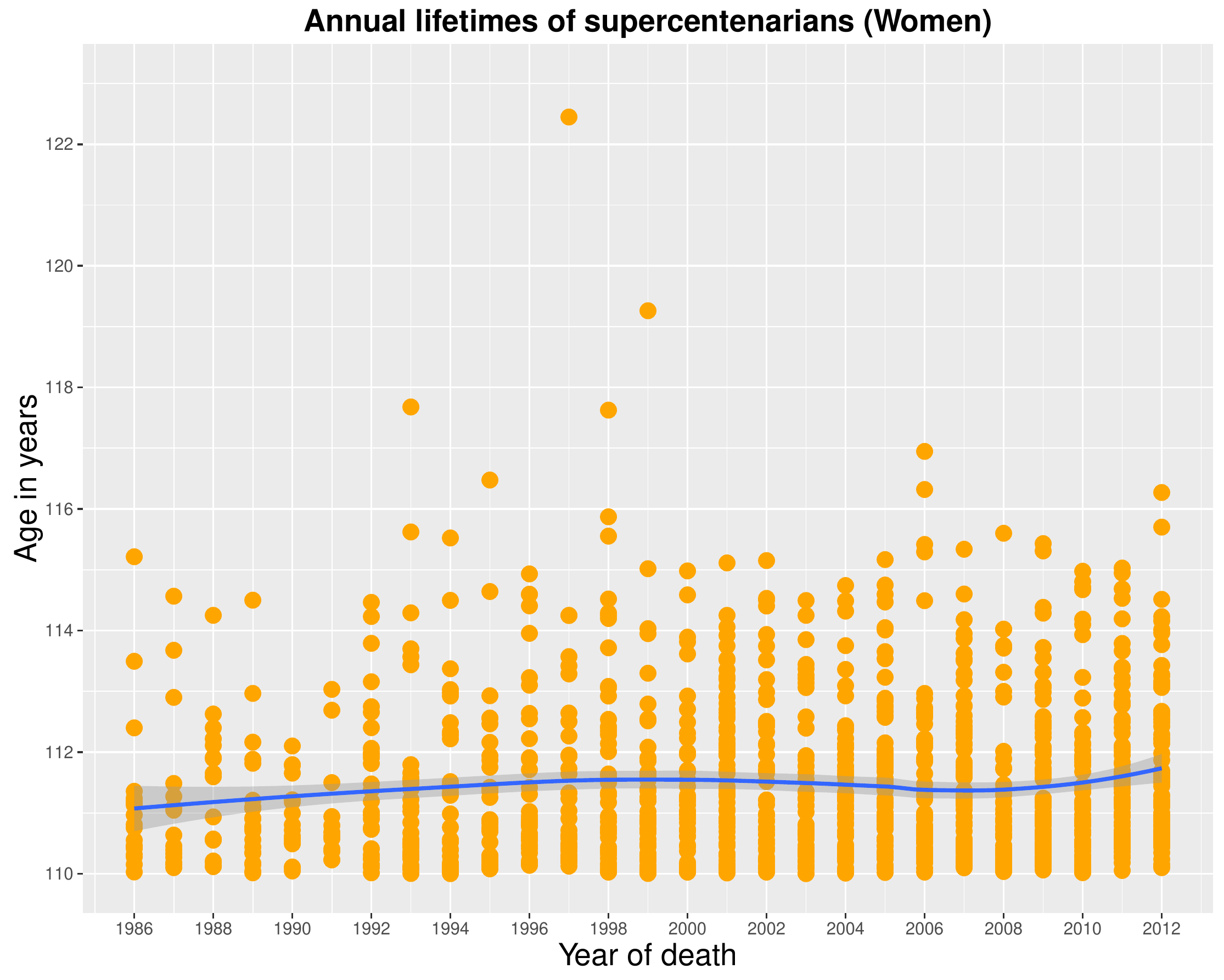}
     \caption{\footnotesize Loess fit to the lifetime of supercentenarians (Women) reported from 1986 to 2012.}
      \label{Loess}
\end{figure}

At this point, several assumptions are needed about the right tail of the lifetime distribution, which is the focus of our extreme value analysis. The first assumption is that the available data comprises a sample of i.i.d. observations.  We find reasonable to assume independence in our data, since we have one record for each individual person. The stationarity assumption is preliminary  asserted from the  plot in Figure \ref{Box_Hist} displaying the comparative boxplots for the larger women's lifetimes by the year. A common feature to all the boxplots is the presence of very large observations classified as extreme outcomes. The boxplots also suggest an increasing third quartile as we progress in time. Such an increase is not so apparent in the annual maxima as we move across the $27$ time points (years). A possible interpretation is that an increase in the mean of the supercentenarians' lifetimes may not be connected to an increasing lifespan over time, but rather to a possible trend in the frequency of the highest lifetime observations. There is a recent semi-parametric development by \citet{dHKTN:15}, suitable for assessing the presence of a trend in the frequency of extreme observations, which is also reflected in the scale of extremes. However, their inference techniques require a large number of replicates per year, and this is not at our grasp given the limited amount of observations available within the same year. For instance in the top part of Figure \ref{Box_Hist}, that we have always less than 40 yearly observations until 1997 (the year of Calment), with the lowest number of 12 observations for 1991. Although the number of observations virtually doubles for the later years, the absolute record still stands on the Calment limit of 122 years, the overall sample maximum. The plot in Figure  \ref{Loess} shows the Loess fit (Local scatterplot smoothing) to the yearly data, given by  blue curve overlay. This is almost parallel to the horizontal axis. Confidence bands are also presented  with a preassigned 95\% confidence level. These seem rather narrow. Hence, we find no evidence of a particular trend through this nonparametric technique.

Despite the above, a POT parametric approach is  applied to detect  a possible trend in the scale. Here, the GPD is fitted  to
the threshold excesses, via a ML fit to $H_\gamma(x/\sigma)$, GPD is fitted to the threshold excesses,  considering  a trend through time as  $\sigma=\sigma_t=\exp\{\beta_0+\beta_1\,t \}$. The resulting non-stationary model is denoted by $\mathcal{M}_1$,  whereas the corresponding stationary  version GPD with $\sigma=\exp\{\beta_0\}$ is denoted by $\mathcal{M}_0$.
Let $\ell_1$  and $\ell_0$ be the maximized log-likelihoods for models  $\mathcal{M}_1$ and $\mathcal{M}_0$, respectively. The Likelihood Ratio  test for  $H_0: \mathcal{M}_1=\mathcal{M}_0$ using the deviance statistic  $D=2\{\ell_1- \ell_0 \}$ to formally compare models $\mathcal{M}_1$ and $\mathcal{M}_0$, returns the results summarized in Table \ref{table2}, for the lifetime data over  the selected thresholds $110,111,112, 113$ and $114$.
The last column contains the $p$-values, related to the different threshold selection. Again, we find  no strong evidence of a linear trend
in the log-scale parameter $\sigma$. Moreover, for each threshold, the EVI and endpoint ML estimates for women lifetimes
  are also listed for the sake of comparison with subsequent results in Section \ref{xF_EVI}.
  \begin{table}
  \caption{\footnotesize Maximum-likelihood parametric GPD($\sigma_t$) and GPD($\sigma$) fit to the lifetime exceedances of supercentenarian women for the time period 1986-2012.}\label{table2}
\begin{center}
 \resizebox{0.9\textwidth}{!}{ 
\begin{tabular}{l|rrcr|c c|c}
  \hline
  $\quad$ $\sigma_t$ & $\ell_i$, $i=0,1$ & $\hat \beta_0$ & $\hat \beta_1$   & $\hat\gamma$ & thresh$/\#$exc &  $\hat x^F_{POT}$ & $p$-value \\
\hline
$ \exp\{\beta_0+\beta_1\,t\}$&  -1739.324 &    0.327 & 0.006&-0.056& 110/1272& & \\
$\exp\{\beta_0\}$  &-1740.238& 0.432  &  -- & -0.061 &  & 135.36 & 0.1763\\
  \hline
$ \exp\{\beta_0+\beta_1\,t\}$&  -867.215 &    0.497 & -0.003&-0.083& 111/639& & \\
$\exp\{\beta_0\}$  & -867.374 & 0.436  &  -- & -0.078 &  & 130.63 & 0.5729\\
  \hline
  $ \exp\{\beta_0+\beta_1\,t\}$&  -424.197 &  0.558 & -0.013 &-0.062& 112/338& & \\
$\exp\{\beta_0\}$  & -425.542 & 0.304  &  -- & -0.045 &  & 142.05 & 0.1010\\
  \hline
  $ \exp\{\beta_0+\beta_1\,t\}$&  -192.652 &  0.567 & -0.019 & -0.044 & 113/164& & \\
$\exp\{\beta_0\}$  & -193.942 & 0.198  &  -- & -0.016 &  & 191.21 & 0.1083 \\
 \hline
   $ \exp\{\beta_0+\beta_1\,t\}$&  -75.904 &  0.150 &  -0.021 & 0.141 & 114/82& & \\
$\exp\{\beta_0\}$  & -76.495 &  -0.237  &  -- & 0.170 &  & -- & 0.2769 \\
  \hline

 % $ \exp\{\beta_0+\beta_1\,t\}$&   -26.647 &  0.708 &  -0.052 & 0.238 & 115/25& & \\
%$\exp\{\beta_0\}$  & -27.421 &  -0.288  &  -- &  0.385 &  & -- & 0.2133 \\
  \hline
\end{tabular}
}
\end{center}
\end{table}

The previous preliminary parametric analysis does not exhaust all the possible choices of parametric models encompassing a trend in extremes. The interest is not in the selection of the most suitable parametric model for extremes, but in being able to ascertain, with a certain degree of confidence, that dropping out of the time covariate does not affect our subsequent analysis under the assumption of stationary supercentenarian women's lifetimes. For recent applications incorporating information over time, we refer the works of \citet{StephTawn:13} and \citet{dHKTN:15}. The latter comprises a comparative analysis with existing methodologies in a similar context.

\begin{figure}
    \centering
   \includegraphics[width=0.8\textwidth]{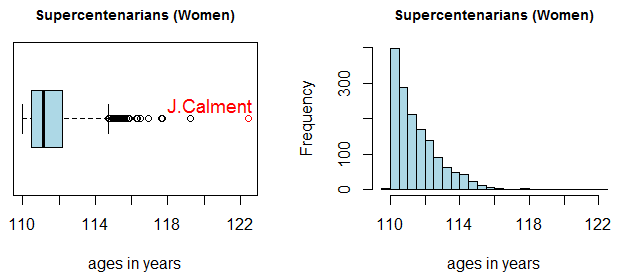}
     \caption{\footnotesize Boxplot and histogram built on the lifetimes of supercentenarian women, for the time interval 1986-2012.}
      \label{Box_Hist0}
\end{figure}

The available supercentenarian women data should be regarded as the greatest lifetimes collection ever attained by the human population. We have $1272$ observations available, which we have found to satisfy the i.i.d. assumption. Figure  \ref{Box_Hist0} contains the boxplot and the histogram for the whole univariate data set.

\subsection{Testing finiteness in the right endpoint}
\label{Testing-xF}

Our first aim is to assess finiteness in the right endpoint of the d.f. $F$ underlying the women lifespan data. The detection of a possibly finite upper bound on our data follows a semiparametric approach, meaning that we essentially assume that $F$ belongs to some max-domain of attraction.  We then consider the usual asymptotic setting, where $k=k_n \rightarrow \infty$ and $k_n/n \rightarrow 0$, as $n\rightarrow \infty$, and hence $X_{n-k,n} \rightarrow x^F\; a.s.$ According to this setting, it is only natural to expect that any statistical approach to the problem of whether there is a finite endpoint or not, will depend on the extent of the dip into the original sample of supercentenarian women's records. The baseline to this issue is mostly driven by a second and more operative question: how to select the adequate top sample fraction to use with both our testing and estimation methods? A suitable choice of $k$ comes from a similar approach to the one in \citet{WANG:95}, where $k_{opt}$ is deemed to be selected at the value $k$ from which  the null hypothesis is rejected.

\begin{figure}
    \centering
   \includegraphics[width=0.65\textwidth]{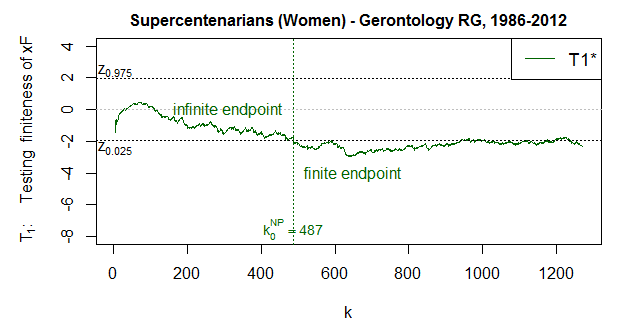}
     \caption{\footnotesize Detecting finiteness in right endpoint for  verified supercentenarians data set: sample paths of the normalized statistics. The horizontal dashed lines correspond to the $\alpha=5\%$ critical barriers.}
      \label{TEST2}
\end{figure}

For a more definite judgment about the existence of a finite upper bound on the supercentenarian lifetimes, we are going to apply the testing procedure introduced in  \citet{ NevesP:10}. The purpose now is to detect finiteness in the right endpoint of the underlying distribution which may belong to either Weibull or Gumbel domains. More formally, the testing problem
 \begin{equation*}
   H_0: F\in{\mathcal D}(G_0),\, x^F=\infty\quad vs \quad H_1: F\in{\mathcal D}(G_\gamma)_{\gamma\leq0},\,x^F<\infty
 \end{equation*}
is tackled using the log-moments $N_r\equiv  N_{n,k}^{(r)},\,\,r=1,2$, defined in \eqref{moms-inv}, but now replacing the observations $X_{i,n}$ by their log-transform $\log(X_{i,n})$. It is possible to do so because we are dealing with positive observations. We point out however that this leads to a non-location invariant method.
The test  statistic $T_1$ being used is defined as
\begin{equation*}
  T_1:=\frac{1}{k}\sum_{i=1}^k\frac{X_{n-i,n}-X_{n-k,n}-T}{X_{n,n}-X_{n-k,n}}, \mbox{ with } T:=X_{n-k,n}\frac{N_1}{2}\left(1-\frac{\left[N_1\right]^2}{N_2}\right)^{-1} \, .
\end{equation*}
 Under $H_0$  the  standardized version of the test, $T_1^*:=\sqrt{k}\,\log k \,T_1$, is asymptotically normal.
Moreover,  $T_1^*$ tends to inflect to the left for bounded tails in the Weibull domain and to the right if the underlying distribution belongs to the Gumbel domain. The rejection region of the test is given by $| T^*_1 |\geq z_{1-\alpha/2}$, for an approximate $\alpha$ significance level.
Figure \ref{TEST2} displays the sample path of $T_1^*$. The most adequate choice of the intermediate number $k$ (which carries over to the subsequent semi-parametric inference) is set on the lowest $k$ at which the critical barriers with a  $\alpha=5\%$ significance level are crossed. This optimality criterium yields $k_{0}^{NP}=487$, spot on the smallest value where we find enough evidence of a finite endpoint.

It remains to be assess whether the distribution underlying the supercentenarian data (now assumed bounded from above) belongs to the Gumbel domain or to the Weibull max-domain of attraction. This will be carried out by a proper hypothesis-testing problem, termed statistical choice of extremes domains. In view of our specific interest on the finite endpoint, we are using the one-sided version of the test. The hypotheses are:
\begin{equation*}
H_0: F\in{\mathcal D}(G_0)\quad vs \quad H_1': F\in{\mathcal D}(G_\gamma)_{\gamma<0}  \,\,.
\end{equation*}

\begin{figure}
    \centering
   \includegraphics[width=0.65\textwidth]{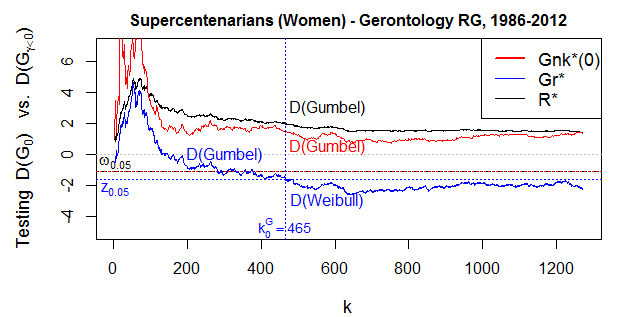}
     \caption{\footnotesize Testing max-domains of attraction with short tailed alternative, for  verified supercentenarians data set:  sample paths of the normalized test statistics. The horizontal dashed lines correspond to the $\alpha=5\%$ critical barriers.}
      \label{TEST1}
      \end{figure}

Figure  \ref{TEST1}  depicts the sample paths  of the new test statistic $G_{n,k}^*(0)$ from Theorem \ref{Thm.Test}, the Greenwood $Gr^*$ statistic,  and the Ratio $R^*$ statistic, as well as  their $\alpha=5\%$ asymptotic critical values. The Greenwood test finds enough evidence to reject the Gumbel domain  hypothesis in favour of a bounded short-tail in Weibull domain, at the suitable intermediate value of $k_{0}^{G}=465$. The two other statistics, new $G_{n,k}^*(0)$ and ratio $R^*$, lead to a more conservative conclusion, with both tests leaning towards the non rejection of the null hypothesis. This conservative aspect also crops up in the simulations section (Section \ref{SecSims}), where the Greenwood test is found to be more powerful against  short-tailed alternatives attached to $\gamma<0$.

From the previous analysis, we find it reasonable to assume a finite right endpoint for some distribution  in the Weibull domain simultaneously for all the adopted testing methods at the maximum number of upper extremes
\begin{equation*}
	 k_{opt}:=\max\{k_{0}^{G}, k_{0}^{NP} \}=487.
\end{equation*}

Therefore, based not only on the testing procedures but also on the complementary EVI estimation presented in Section \ref{xF_EVI}, it seems reasonable to conclude that the lifespan distribution belongs to the Weibull domain of attraction with finite right endpoint $x^F$.\\

For definitely discarding the presence of a heavy-tailed distribution underlying our data set, it is also important to test:
 \begin{equation*}
H_0: F\in{\mathcal D}(G_0)\quad vs \quad H_1: F\in{\mathcal D}(G_\gamma)_{\gamma>0}  \,.
\end{equation*}

\begin{figure}
    \centering
   \includegraphics[width=0.65\textwidth]{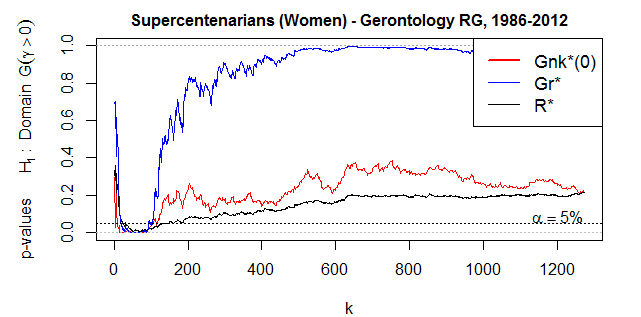}
     \caption{\footnotesize Testing max-domains of attraction for  verified supercentenarians data set: $ p$-values   for   the normalized test statistics, with heavy tailed alternatives. The horizontal dashed line corresponds to the $\alpha=5\%$ nominal level.}
      \label{TEST10}
      \end{figure}

In Figure \ref{TEST10} we observe that all the $p$-values determined by all the three test statistics are increasing with $k$. Again, the conservative behavior of the two tests seems to emerge. Despite all $p$-value paths begin at very small values around zero, this only lingers for a tight range of higher thresholds which may not be in good agreement with the requirement of a sufficiently large $k$. The Greenwood statistic returns $p$-values very close to 1, from about $k=400$ onwards. The other two statistics ($G_{n,k}^*(0)$ and $R^*$) yield moderate  $p$-values, with larger values returned by $G_{n,k}^*(0)$. It seems sensible then to discard a heavy-tailed distribution for the supercentenarian women lifespan, a conclusion clearly verified by the Greenwood test.

\subsection{Endpoint estimation for women's records without EVI knowledge} \label{xF}

Following the testing procedures in \ref{Testing-xF} and the reported optimal number $k_{opt}=487$, we will present analogous graphical tools with respect to endpoint estimation. The first purpose is to illustrate the  smooth behaviour of $\hat{x}^F$ defined in \eqref{EPEst}, already anticipated in the simulations section (Section \ref{SecSims}).

\begin{figure}
    \centering
   \includegraphics[width=0.7\textwidth]{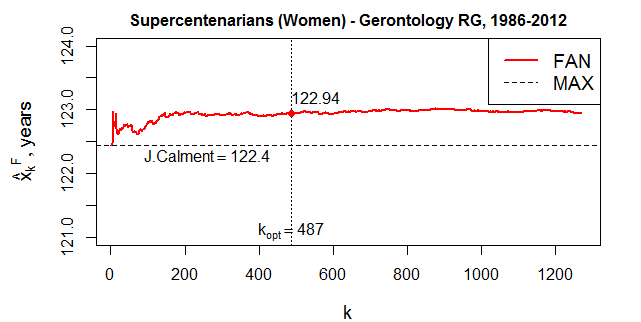}
     \caption{\footnotesize Endpoint estimation for  verified supercentenarians data set: estimator $\hat x^F$ and Calment limit.}
      \label{xF_data0}
\end{figure}

Figure \ref{xF_data0} displays the comparative finite-sample behaviour of $\hat{x}^F$ (notation: FAN) with the na\"ive Calment limit (notation: MAX) for the supercentenarian's data set. Recall that our  database is regarded as a collection of the greatest lifetimes of women population and that the general endpoint estimator $\hat{x}^F$ always returns values above the na\"ive endpoint estimate, i.e. greater than ``Calment limit" of $122.4$ years.

After the initial rough path in the range of approximately one hundred top observations, the estimates trajectory of $\hat{x}^F$ then becomes very flat. Once we dip into the intermediate range of extremes, the sample path becomes smoother.

At this point we find it sensible to provide similar information provided by the other semi-parametric and parametric endpoint estimators already intervening in the simulation study, as they constitute a good complement to a more thorough insight about the true endpoint. This is the subject of the next section.

\subsection{Linking endpoint estimation of women's records to EVI estimation} \label{xF_EVI}
 \label{RB2_CI}

\begin{figure}
    \centering
  \includegraphics[width=0.7\textwidth]{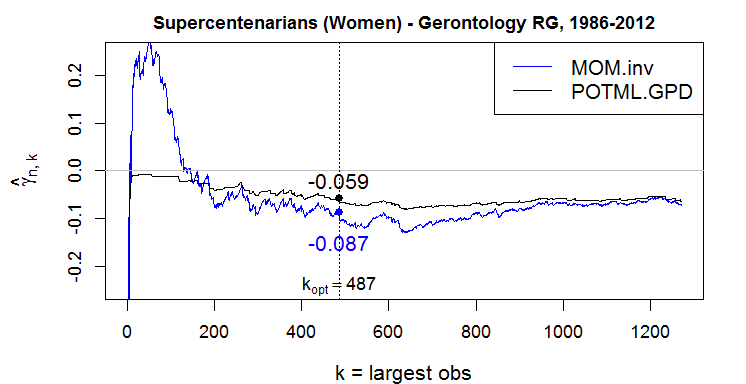}
     \caption{\footnotesize  EVI estimation with the verified supercentenarian data.}
      \label{EVI_data}
\end{figure}

In this section, the endpoint estimation is tackled using methodologies that require an estimator for  the extreme value index $\gamma<0$. In this sequence, we adopt the moment related estimator $\hat{\gamma}_{n,k}^{-}$ (notation: MOM.inv), defined in \eqref{EVI-inv}, and the POTML.GPD estimator of the shape parameter $\hat{\gamma}_{n,k}^{ML}$ to plug in \eqref{POTML} (subject to $\gamma<0$). Figure \ref{EVI_data} is the estimates plot of the chosen estimators for $\gamma$, as function of $k$.  We note that these two estimators enjoy the location and scale invariance property. Retaining the $k_{opt}=487$ larger observations (a value delivered by the testing procedures) as the effective sample, we find the point estimates $\hat{\gamma}_{n,487}^{-}=-0.087$  and $\hat{\gamma}_{n,487}^{ML}=-0.059$, both coherent with a short-tailed distribution attached to some  $\gamma \in (-1/2,0)$.

\begin{figure}
    \centering
    \includegraphics[width=0.7\textwidth]{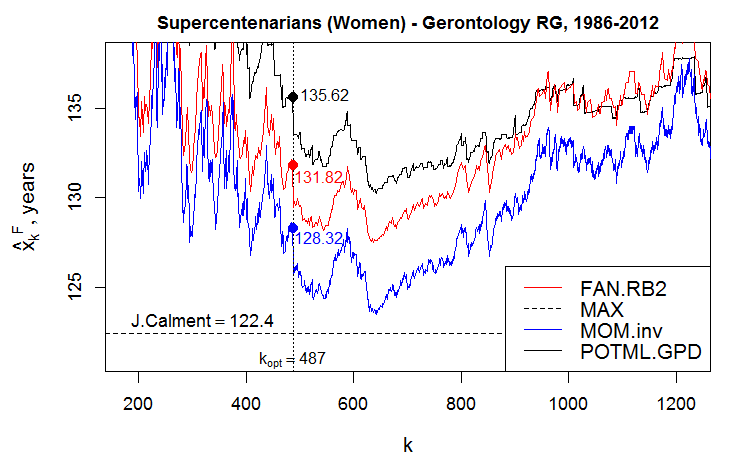}
     \caption{\footnotesize Endpoint estimation with the verified supercentenarian data.}
      \label{xF_data1}
\end{figure}

Figure \ref{xF_data1} depicts the results for estimators $\hat x^F_{ML}$, $\hat x^*$ (defined in \eqref{POTML} and \eqref{mom_est}, respectively). We are also including the second order reduced bias version of the general endpoint estimator $\hat x^F_{RB2}$ defined in \eqref{RB2}, by plugging in the estimators $\hat\gamma=\hat\gamma^{-}_{n,k}$ and $\hat a_0\bigl(\ndivk\bigr)$ provided in \eqref{EVI-inv} and \eqref{ank-inv}, respectively. The  trajectory of $\hat x^F_{RB2}$ lies in the  middle range between that of $\hat x^F_{ML}$ and of $\hat x^*$.

The optimal value $k_{opt}=487$ can viewed as benchmark value (or change point) since it breaks the disruptive estimates pattern. It actually pinpoints where the graph stops being too rough to make inference and starts being more stable,  so that we can  infer about the endpoint. The latter applies in particular to estimators $\hat x^*$,  $\hat x^F_{RB2}$ and  $\hat x^F_{ML}$, which return $\hat x^*_{k_{opt}}=128.32$,  $\hat x^F_{RB2,k_{opt}}=131.82$ and  $\hat x^F_{ML,k_{opt}}=135.62$. The simulations also outline the plain general endpoint estimator $\hat x^F$ as being relatively efficient (in terms of bias and MSE) on those moderate values of $k^*$ where other estimators fall short. In this respect, Figure \ref{xF_data0} shows a rather flat plateau from which we find safe to draw an estimate for the endpoint, and this portion of the graph includes $\hat x^F_{k_{opt}}=122.94$. The asymptotic results in Theorem \ref{MainCor} for $\gamma \in (-1/2,0)$ equip us  with a tool for finding an approximate $100(1-\alpha)\%$ upper bound for the true endpoint $x^F$. Obviously, this process calls for the estimation of $\gamma$ since we need some guidance about the proper interval where the EVI lies within. Hence, we have just introduced the most direct link to the estimation of $\gamma$ (see Figure \ref{EVI_data}) in connection with the general endpoint estimator $\hat{x}^F$. Selecting $k_{opt}=487$, the $95\%$ confidence upper bound for $x^F$ delivered by \eqref{up_bound}  is $133.23$ years.\\

We also note that the simulation outcomes relate well to the present results arising from the available data of supercentenarian women's records. For instance, we observe in Figure \ref{EVI_data} that small values of $k$ ($k\leq 110$, approximately) find positive estimates for the MOM.inv estimator which could on their own account for the erratic pattern of $\hat x^*_{k}$ in the plot of Figure \ref{xF_data1}, but the simulations have yield this rough pattern very often in connection with a true negative EVI. Furthermore, the shape parameter $\gamma$ is estimated via POTML subject to $\gamma <0$ (cf. Figure \ref{EVI_data}) and this returns endpoint estimates $\hat x^F_{ML}$ well aligned with the ones pertaining to  $\hat x^F_{RB2}$ and $\hat x^*_{k}$ (cf. Figure \ref{xF_data1}). Therefore, heeding to the simulation results, we find reasonable to conclude that the misrepresentation of the negative EVI for moderate values of $k$  is not the main factor compromising  the performance of the adopted endpoint estimators in the early part of the plot (amounting to about $10\%$ of the original sample size $n$, say).

Finally, this practical application also seems to suggest that removing the bias component from $\hat x^F$ causes an increase in the variance. The bias/variance trade-off effect is grasped more thoroughly in the Appendix \ref{RB2_MC}, where the finite sample properties of the second order reduced bias estimator $\hat x^F_{RB2}$ are studied by taking Models 1 and 4, from Section \ref{SecSims}, as  parent distributions. The brief simulation study in the Appendix \ref{RB2_MC} is expected to reinforce the suggested competitive performance amongst the above-mentioned endpoint estimators.

%Estimated  $95 \%$ - one-sided CI   133.23\\
%Estimated  $99 \%$ - one-sided CI   133.64\\
%Estimated  $99.5 \%$ - one-sided CI   133.78\\

% ============================================================
\subsection{An upper limit to \emph{lifespan} and probability of surpassing Calment limit}\label{uplimit}

From  the previous data analysis, one would say that the ultimate human lifespan would not be greater than 133.23 years (the estimated upper bound from \eqref{up_bound}, obtained in section \ref{RB2_CI}). This gives some insight beyond  Calment's achievement: the absolute record of 122.4 years, still holding to the present date (for the last 19 years). We thus expect that the probability of exceeding the ``Calment limit", even for a supercentenarian women, will be extremely low.

\begin{figure}[!htbp]
    \centering
  \includegraphics[width=0.6\textwidth]{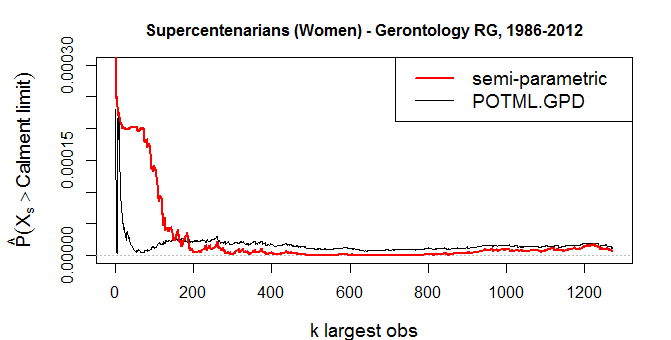}
     \caption{\footnotesize Probability of exceeding the ``Calment limit",  for a supercentenarian women, given today's state of the art.}
      \label{prob}
\end{figure}
\begin{figure}[!htbp]
    \centering
   \includegraphics[width=0.6\textwidth]{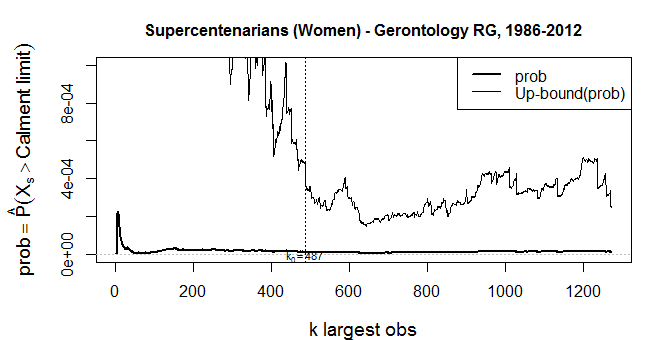}
  \includegraphics[width=0.6\textwidth]{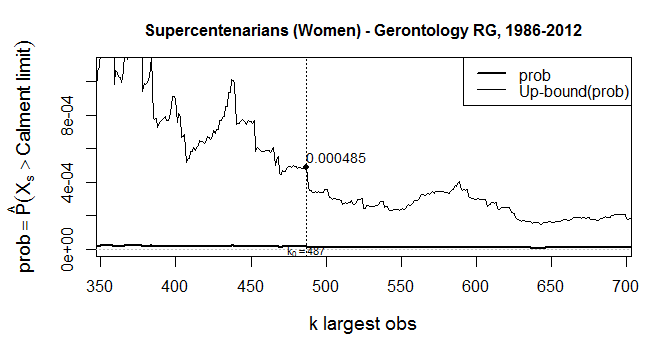}
     \caption{\footnotesize Probability of exceeding the ``Calment limit" using the POTML.GPD approach and corresponding $95\%$ confidence bound: all sample (\emph{top}); $350\leq k\leq 700$ (\emph{bottom}).}
      \label{boundML}
\end{figure}

This tail probability can be estimated using the following semi-parametric estimator:
\begin{equation}\label{prob_tail}
 \widehat P_n (X_{s}>122.4 ) := {k \over n} \left\{\max\left(0,1+\hat{\gamma}_{n,k}^{-} { 122.4- X_{n-k,n} \over
 \hat{a}(n/k) }\right)  \right\}^{-1/\hat{\gamma}_{n,k}^{-}},
\end{equation}
\citep[cf. (4.4.1) in][]{deHaanFerreira:06} where $X_s:= X|X \geq 110$ denotes the lifetime of a supercentenarian women, and $\hat{\gamma}_{n,k}^{-}$, $  \hat{a}(n/k) $ are the related estimators defined in \eqref{EVI-inv} and \eqref{ank-inv}, respectively.

Figure \ref{prob} depicts the probability estimates from \eqref{prob_tail}, together with their POTML.GPD analogues, for a wide range of larger values of  the 1272 verified supercentenarians  data set. In contrast with the previous statistical analysis,   the sample size $n$ now intervenes in \eqref{prob_tail}. Therefore, any inference drawn on this account will apply to the subpopulation of  supercentenarian women under study. All point estimates are very close to zero. Figure \ref{boundML} displays again the POTML.GPD estimates but with respective $95\%$ upper bounds. Since we are dealing with very small probabilities the asymptotic bounds are not so sharp, as opposed to the case of an underlying distribution with infinite right endpoint.

\section{Concluding remarks} %========== Sec. Conclusions =================
\label{SecConclusion}

The scope for application of the right endpoint estimator introduced in \citet{FAN:14}, primarily designed for the Gumbel domain of attraction, is here extended to the case of an underlying distribution function $F$ in the Weibull domain. The consistency property and asymptotic distribution of this general endpoint estimator $\hat{x}^F$ renders a unified estimation procedure for the right endpoint under the assumption that $F\in \mathcal{D}(G_{\gamma})_{\gamma\leq 0}$. A new test statistic arises tied-up with $\hat{x}^F$ thus  incrementing the range of available testing procedures for selecting max-domains of attraction.
Our main findings are listed below.

\begin{itemize}

\item The general endpoint estimator does not require the estimation of the  EVI, unlike the widely-used semi-parametric alternatives.
    \item By construction, the estimator $\hat{x}^F$ always  returns larger values than the sample maximum $X_{n,n}$, a property not shared by other semi-parametric methodologies we have encountered so far, particularly those predicated on the Weibull max-domain of attraction.
    \item  The simulation study conveys a good finite sample performance of  the general endpoint estimator, ascertaining competitiveness to benchmark endpoint estimators specifically tailored for the Weibull domain.
    \item Related to the previous, the general endpoint estimator performs better for distributions with some $\gamma>-1/2$, which corresponds to the most common situation in practical applications.

 \item The problem of choosing the most adequate number $k$ of upper order statistics is here  mitigated by the usual flat pattern of the estimates trajectories, a typical feature of the general endpoint estimator.

 \item The application to the supercentenarian women's lifetimes illustrates how we can easily establish a confidence upper bound to the right endpoint, building on the asymptotic results for $\gamma>-1/2$.
 \end{itemize}

\section{APPENDIX A: Proofs} %============= APPENDIX: PROOFS =================================
\label{SecProofs}

This section is entirely dedicated to the proofs of the results introduced in Section \ref{SecEstim}.
In what follows we find more convenient to consider the estimator $\hat{x}^F$ in the functional form
\begin{equation}\label{Funct}
	\hat{x}^F= X_{n-k,n}-\frac{1}{\log 2}\intunit \bigl( X_{n-[2ks],n}-X_{n-[ks],n}\bigr)\, \frac{ds}{s},
\end{equation}
where $[a]$ denotes the integer part of $a\in \real$ \citep[more details about the representation \eqref{Funct} can be obtained in][]{FAN:14}.

We note that if  $s\in [0,\,1/(2k)[$, then the integral in \eqref{Funct} is equal to zero. Bearing this in mind, we write
\begin{equation*}
\hat{x}^F= X_{n-k,n}-\frac{1}{\log 2}\intab{\frac{1}{2k}}{1} \bigl( X_{n-[2ks],n}-X_{n-[ks],n}\bigr)\, \frac{ds}{s}.
\end{equation*}
Moreover, if $s\in [1/(2k),1/k[$ then $[ks]=0$ (not depending on $s$)  and thus $X_{n-[ks],n}=X_{n,n}$. Therefore, we have that
\begin{equation}\label{FunctAux}
\hat{x}^F= X_{n-k,n}-\frac{1}{\log 2} \Bigl\{\intab{\frac{1}{2k}}{\frac{1}{k}} X_{n-[2ks],n}\, \frac{ds}{s}-X_{n,n}\,\intab{\frac{1}{2k}}{\frac{1}{k}} \frac{ds}{s}+\intab{\frac{1}{k}}{1}( X_{n-[2ks],n}-X_{n-[ks],n})\,\frac{ds}{s}\Bigr\}.
\end{equation}
With a suitable variable transform  on the  last integral, we can reassemble \eqref{FunctAux} in a tidy manner:
\begin{equation}\label{FunctClean}
	\hat{x}^F = X_{n,n}+X_{n-k,n}-\frac{1}{\log 2}\intab{\frac{1}{2}}{1} X_{n-[2ks],n}\, \frac{ds}{s}.
\end{equation}
This is the main algebraic expression that will be used to derive the asymptotic distribution of $\hat{x}^F$
in the proof of Theorem \ref{ThmMain}, which is a natural consequence of the three random contributions in \eqref{FunctClean}.\\\\

%\begin{proof}
\noindent
\textbf{Proof of Proposition \ref{PropCons}:}
 We see that the integral in the functional form \eqref{FunctClean} satisfies the inequalities
\begin{equation*}
	(\log 2 )\,  X_{n-k,n} \leq \intab{\frac{1}{2}}{1} X_{n-[2ks],n}\, \frac{ds}{s} \leq ( \log 2) \,  X_{n-2k,n}.
\end{equation*}
 Therefore, we obtain the following upper and lower bounds involving $\hat{x}^F-x^F$,
\begin{equation*}
X_{n,n}-x^F\, \leq \, \hat{x}^F - x^F \, \leq \,  (X_{n,n}-x^F)+ X_{n-k,n}-X_{n-2k,n},
\end{equation*}
and the result thus follows easily because the three o.s. $X_{n,n}$, $X_{n-k,n}$ and $X_{n-2k,n}$ all converge almost surely to $x^F$, provided the intermediate nature of $k=k_n$.
%\end{proof}
\qed

      \begin{rem}
Alternative proof  based on the functional form \eqref{EPEst_alter} of the $k^*:=2k$ top o.s.: strong consistency of the general endpoint estimator comes from the lower and upper bounds of \eqref{EPEst_alter} given below.
  \begin{eqnarray*}
    \hat{x}^F - x^F &=& (X_{n,n} - x^F) + \left( X_{n-k,n} -  \frac{1}{\log 2} \sum_{i=0}^{k-1} \log(\frac{k+i+1}{k+i})\, X_{n-k-i,n}\right) \\
  & \geq & (X_{n,n} - x^F) + \left( X_{n-k,n} -  X_{n-k,n}\, \frac{1}{\log 2} \sum_{i=0}^{k-1}\log(\frac{k+i+1}{k+i})\right)= X_{n,n} - x^F
  \end{eqnarray*}
and on the other hand,
 \begin{eqnarray*}
\hat{x}^F - x^F   & \leq & (X_{n,n} - x^F) + \left( X_{n-k,n} - X_{n-2k+1,n}\, \frac{1}{\log 2}\sum_{i=0}^{k-1}\,\log(\frac{k+i+1}{k+i}) \right)\\
            & = & (X_{n,n} - x^F) +( X_{n-k,n} -  X_{n-2k+1,n}) \, \,       ;
\end{eqnarray*}
since for any intermediate $k=k_n$ the o.s. $X_{n,n}, X_{n-k,n},X_{n-2k,n} $ converge
almost surely to $x^F$, the result follows.
\end{rem}

Before getting under way to the proof of the main Theorem, we need to lay down some ground results. These comprise  a Proposition regarding $\gamma <0$ and a Lemma for general $\gamma$.

%--------- Proposition ----------------------------------

\begin{prop}\label{PropxF}
Suppose $X_{n,n}$ is the maximum of a random sample whose parent d.f. $F$ detains finite right endpoint of $F$, i.e. $x^F=U(\infty) < \infty$. Assume the second order condition \eqref{2ERVU} holds with $\gamma <0$. If $k=k_n$ is such that, as $n\rightarrow \infty$, $k \rightarrow \infty$, $k/n \rightarrow 0$ and $\sqrt{k}\,A_0(n/k) \rightarrow \lambda^* \in \real$,  then
\begin{enumerate}
\item for $\gamma \geq -1/2$, for each $\varepsilon >0$,
	\begin{equation}\label{Aux1}
k^{-\gamma-\varepsilon}\biggl|\frac{X_{n,n}-x^F}{a_{0}\bigl(\ndivk\bigr)}\biggr| \,\conv{p} 0.
\end{equation}
Moreover,
\begin{equation*}
	k^{-\gamma} \, \frac{X_{n,n}-x^{F}}{a_0\bigl(\ndivk\bigr)}
	\, \conv{d}\,  \frac{Z^{\gamma}}{\gamma},
	\end{equation*}
	where  $Z$ denotes a standard Fr\'echet with d.f. $\Phi_{1}$ as in \eqref{Frech}.
\item for $\gamma < -1/2$,
\begin{equation*}
\sqrt{k}\,\biggl|\frac{X_{n,n}-x^F}{a_{0}\bigl(\ndivk\bigr)}\biggr| \,\conv{p} 0.
\end{equation*}
\end{enumerate}
\end{prop}

%\begin{proof}
\noindent
\textbf{Proof:} Owing to the well-known equality in distribution that $X_{i,n}\id \,U(Y_{i,n})$, $i=1,2,\ldots,n$, with $\bigl\{Y_{i,n}\bigr\}_{i=1}^n$ the $n$-th o.s. from a sample of $n$ independent r.v.s with common (standard) Pareto d.f. given by $1-x^{-1}$, $x\geq 1$, then the following equality in distribution holds:
\begin{equation*}
	\frac{X_{n,n}-x^{F}}{a_0\bigl(\ndivk\bigr)}
	\, \id \,  \Bigl\{\frac{U\bigl(\kdivn Y_{n,n}\,\ndivk\bigr)-U\bigl(\ndivk\bigr)}{a_0\bigl(\ndivk\bigr)}+
\frac{1}{\gamma}\Bigr\}-\Bigl\{\frac{U(\infty)-U\bigl(\ndivk\bigr)}{a_0\bigl(\ndivk\bigr)}+\frac{1}{\gamma}\Bigr\}.
	\end{equation*}
Now we use conditions \eqref{2ERVU}  and \eqref{2ERVinfty} with $t$ replaced by $n/k$ everywhere:
\begin{eqnarray*}
	\frac{X_{n,n}-x^{F}}{a_0\bigl(\ndivk\bigr)} &\id& \biggl\{\frac{k^{\gamma}\bigl(n^{-1}Y_{n,n}\bigr)^{\gamma}-1}{\gamma}+\frac{1}{\gamma}+A_0\bigl(\ndivk\bigr)\Psi^{\star}_{\gamma,\rho}\Bigl(\kdivn \,Y_{n,n}\Bigr)\bigl(1+o_p(1)\bigr)\biggr\}\\
	& & \mbox{\hspace{4.0cm}} -\biggl\{A_0\bigl(\ndivk\bigr)\Psi^{\star}_{\gamma,\rho}(\infty)\bigl(1+o(1)\bigr)\biggr\}\\
	&=& \frac{k^{\gamma}\bigl(n^{-1}Y_{n,n}\bigr)^{\gamma}}{\gamma}+ A_0\bigl(\ndivk\bigr)\Bigl\{\Psi^{\star}_{\gamma,\rho}\Bigl(\kdivn \,Y_{n,n}\Bigr)+\frac{1}{\gamma+\rho} I_{\{\rho<0\}}\Bigr\}+o_p\Bigl(A_0\bigl(\ndivk\bigr)\Bigr)
\end{eqnarray*}
We note at this stage that  $n^{-1}Y_{n,n}$ is asymptotically a Fr\'echet r.v.  with d.f. given by $\Phi_{1}$ in \eqref{Frech}. This non-degenerate limit yields $(k/n)Y_{n,n}$ going to infinity with probability one, which implies in turn that
$\Psi^{\star}_{\gamma,\rho}\Bigl(k \bigl(n^{-1}Y_{n,n}\bigr)\Bigr)\rightarrow -\,(\gamma+\rho)^{-1} I_{\{\rho<0\}}$, as $n\rightarrow \infty$. Therefore, we obtain for $\gamma \geq -1/2$,
\begin{equation}\label{Max}
k^{-\gamma}\frac{X_{n,n}-x^{F}}{a_0\bigl(\ndivk\bigr)}\,\id \, \frac{\bigl(n^{-1}Y_{n,n}\bigr)^{\gamma}}{\gamma}+ o_p\bigl(k^{-\gamma-1/2}\bigr),
\end{equation}
by virtue of $\sqrt{k}A_0(n/k)=O(1)$, and \eqref{Aux1} thus follows directly for each $\varepsilon >0$. The second part in point 1. is ensured from \eqref{Max} by the continuos mapping theorem. For $\gamma < -1/2$, we observe from \eqref{Max} that
\begin{equation*}
\sqrt{k}\,\frac{X_{n,n}-x^{F}}{a_0\bigl(\ndivk\bigr)}\,\id \, k^{1/2+\gamma}\frac{\bigl(n^{-1}Y_{n,n}\bigr)^{\gamma}}{\gamma}+ o_p(1).
\end{equation*}
Since we are addressing the case $\gamma+1/2<0$, the fact that $n^{-1}Y_{n,n}$ converges in distribution to a Fr\'echet r.v. suffices to conclude the proof.
%\end{proof}
\qed  % flush right qed marks, e.g. at end of proof
 %------------------------------------------------------------

 \begin{lem}\label{LemVector}
Suppose that $U$ satisfies the second order condition \eqref{2ERVU} with $\gamma \in \real$ and $\rho \leq 0$. If $k=k_n$ is an intermediate sequence such that $\sqrt{k}\, A_0(n/k) =O(1)$, then
\begin{equation}\label{RandVect}
	\sqrt{k}\,\bigl(P_n,\, Q_n\bigr):= \sqrt{k}\,\biggl( \intab{1/2}{1}\frac{X_{n-[2ks],n}-U\bigl(\frac{n}{2ks}\bigr) }{a_0\bigl(\ndivk\bigr)}\, \frac{ds}{s},\, \frac{X_{n-k,n}-U\bigl(\ndivk\bigr)}{a_0\bigl(\ndivk \bigr)}\biggr)
\end{equation}
converges in distribution to the bivariate normal $(P,\,Q)$ random vector with zero mean
and covariance structure given by
\begin{eqnarray*}
E(P^2)&=& \left\{ \begin{array}{ll}
                                \frac{2}{\gamma}\Bigl(\frac{2^{-(2\gamma+1)}-1}{2\gamma+1}-\frac{2^{-(\gamma+1)}-1}{\gamma+1}\Bigr) , & \mbox{ }  \gamma \neq 0, \\
                                1-\log 2, & \mbox{ }  \gamma= 0,
			 \end{array}
                            \right.\\
E(P\,Q) &=& \left\{ \begin{array}{ll}
                                -\frac{1}{\sqrt{2}}\, \frac{2^{-\gamma}-1}{\gamma}, & \mbox{ }  \gamma \neq 0, \\
                               \frac{\log 2}{\sqrt{2}}, & \mbox{ }  \gamma= 0,
			 \end{array}
                            \right.\\
E(Q^2) &=&  1.
\end{eqnarray*}
\end{lem}

%\begin{proof}
\noindent
\textbf{Proof:}
The first component in \eqref{RandVect} shall be tackled by Theorem 2.4.2 of \citet{deHaanFerreira:06} with $k$ replaced by $2k$ therein. In particular,
\begin{eqnarray}\label{Aux}
\nonumber	 & & 	\sqrt{2k}\,\intab{1/2}{1}\frac{X_{n-[2ks],n}-U\bigl(\frac{n}{2ks}\bigr) }{a_0\bigl(\ndivk\bigr)}\, \frac{ds}{s}\\
	 & =&  \frac{a_0\bigl(\frac{n}{2k}\bigr)}{a_0\bigl(\frac{n}{k}\bigr)}\, \sqrt{2k} \intab{1/2}{1} \Bigl\{ \frac{X_{n-[2ks],n}-U\bigl(\frac{n}{2k}\bigr)}{a_0\bigl(\frac{n}{2k}\bigr)}-\frac{U\bigl(\frac{n}{2ks}\bigr)-U\bigl(\frac{n}{2k}\bigr)}{a_0\bigl(\frac{n}{2k}\bigr)} \Bigr\}\,\frac{ds}{s}
\end{eqnarray}
Then, under the second order conditions \eqref{2ERVU} and \eqref{2RVa}, Theorem 2.4.2 of \citet{deHaanFerreira:06}  yields for the definite integral on the right hand-side of \eqref{Aux}:
\begin{eqnarray*}
	 & & \sqrt{2k}\,\intab{1/2}{1}\frac{X_{n-[2ks],n}-U\bigl(\frac{n}{2ks}\bigr) }{a_0\bigl(\ndivk\bigr)}\, \frac{ds}{s}\\
	 &=& \frac{1}{2^{\gamma}} \intab{1/2}{1} \Bigl\{s^{-\gamma-1}W_n(s) +o_p(1)s^{-\gamma-1/2-\varepsilon}+ o\Bigl(\sqrt{2k}\,A_0\bigl(\frac{n}{2k}\bigr)\Bigr) \Bigr\}\,\frac{ds}{s} +O_p\Bigl(A_0\bigl(\frac{n}{k}\bigr)\Bigr),
\end{eqnarray*}
where $\{W_n(s)\}_{n\geq 1}$, $s>0$, denotes a sequence of Brownian motions.
Under the assumption that $\sqrt{k}A_0\bigl(n/(2k)\bigr)= O(1)$, we  obtain as $n \rightarrow \infty$,
\begin{equation*}
	\sqrt{k}\,P_n=\frac{1}{\sqrt{2}}\intab{1/2}{1} (2s)^{-\gamma }W_n(s)\,\frac{ds}{s^2} +O_p\Bigl(A_0\bigl(\frac{n}{k}\bigr)\Bigr)+o_p(1).
\end{equation*}
If $\gamma =0$, the integral on the right hand side becomes $\int_{1/2}^1 W_n(s)\, ds/s^2$.
In either case, this integral corresponds to the sum of asymptotically multivariate normal random variables.
Now, the second component of the random vector \eqref{RandVect} is asymptotically standard normal \citep[cf. Theorem 2.4.1 of][]{deHaanFerreira:06}.
Finally, the covariance for the limiting bivariate normal, $E(P\,Q)$, is calculated in a straightforward way using similar calculations to the ones in p.163 of \citet{deHaanFerreira:06}.
%\end{proof}
\qed  % flush right qed marks, e.g. at end of proof

%--------- Proof Main Theorem ----------------------------

%\begin{proof}
\noindent
\textbf{Proof of Theorem \ref{ThmMain}:}
Let $h(\gamma)= (\log 2)^{-1} \int_{1/2}^1\bigl\{(2s)^{-\gamma}-1)/(-\gamma)\bigr\}\,ds/s$, which  is defined in \eqref{aga}. Taking the auxiliary function $a_0$ from the second order condition \eqref{2ERVU} we write the following normalization of $\hat{x}^F$ (cf. \eqref{FunctClean} and \eqref{Aux}):
\begin{equation*}
\frac{\hat{x}^F-x^F}{a_0\bigl(\ndivk\bigr)}-h(\gamma)= W_n - \frac{1}{\log2} P_n+Q_n-\frac{1}{\log 2}\intab{1/2}{1}\left(\frac{U\bigl(\frac{n}{2ks}\bigr)-U\bigl(\ndivk\bigr)}{a_0\bigl(\ndivk\bigr)}
-\frac{(2s)^{-\gamma}-1}{\gamma}\right)\,\frac{ds}{s},
\end{equation*}
with $(P_n, Q_n)$ defined in Lemma \ref{LemVector} and $W_n :=\bigl(X_{n,n} - x^F \bigr)/a_0(n/k)$.
Now, Lemma \ref{LemVector} entails that  $\sqrt{k}(P_n,\, Q_n)$ is asymptotically bivariate normal distributed as $(P,Q)$. Proposition \ref{PropxF} expounds the limiting distribution of $W_n$ provided suitable normalization, possibly different than $\sqrt{k}$. Hence, the crux of the proof is in the following distributional expansion, under the second order condition \eqref{2ERVU}, for large enough $n$:
\begin{equation}\label{crux}
k^{-\gamma}\, \Bigl(\frac{\hat{x}^F-x^F}{a_0\bigl(\ndivk\bigr)}-h(\gamma)\Bigr)
	=  k^{-\gamma}W_n  + k^{-(\gamma+1/2)} \Bigl\{ \sqrt{k}Q_n- \frac{ \sqrt{k}}{\log2} \Bigl(P_n+A_0\bigl(\ndivk \bigr)\intab{1/2}{1}\Psi^{\star}_{\gamma,\rho} \bigl(\frac{1}{2s}\bigr)\,\frac{ds}{s} \Bigr)\Bigr\}.
\end{equation}
We shall consider the cases $\gamma >-1/2$, $\gamma =-1/2$ and $\gamma <-1/2$ separately.
\begin{description}
\item[Case $\gamma >-1/2$:] Proposition \ref{PropxF}(1) and Lemma \ref{LemVector} upon \eqref{crux} ascertain the result, by virtue that $W= Z^{\gamma}/\gamma$ with $Z$ a standard Fr\'echet random variable.
\item[Case $\gamma =-1/2$:] The random  component $W_n$ is asymptotically independent of the remainder $P_n$ and $Q_n$. This claim is supported on Lemma 21.19 of \citet{vanderVaart:98}. Hence, the combination of Proposition \ref{PropxF} with Lemma \ref{LemVector} ascertains the result.
\item[Case $\gamma <-1/2$:] It is  now convenient to rephrase \eqref{crux} with a suitable normalization in view of Proposition \ref{PropxF} and the precise statement thus follows:
	\begin{equation*}
		 \sqrt{k}\, \Bigl(\frac{\hat{x}^F-x^F}{a_0\bigl(\ndivk\bigr)}-h(\gamma)\Bigr)
	=   \sqrt{k}\Bigl\{Q_n- \frac{ P_n}{\log2} -\,A_0\bigl(\ndivk \bigr) \frac{1}{\log 2}\intab{1/2}{1}\Psi^{\star}_{\gamma,\rho} \bigl(\frac{1}{2s}\bigr)\,\frac{ds}{s}\Bigr\} +O_p(k^{\gamma+1/2}).
	\end{equation*}
\end{description}
%\end{proof}
\qed  % flush right qed marks, e.g. at end of proof
%---------------------------------------------------------------------

%\vspace{0.5cm}
%\begin{proof}
\noindent
\textbf{Proof of Corollary \ref{MainCor}:}
The result follows immediately from Theorem \ref{ThmMain}, provided $W$ and $N$ are independent random variables. Lemma 21.19 of \citet{vanderVaart:98} ensures the latter.
%\end{proof}
\qed  % flush right qed marks, e.g. at end of proof

%------------------------  TEST  ---------------------------------------------

%\vspace{0.5cm}
%\begin{proof}
\noindent
\textbf{Proof of Theorem \ref{Thm.Test}:}
The test statistic
\begin{equation*}
G_{n,k}:={  \hat x^F- X_{n-k,n}   \over  X_{n-k,n}-X_{n-2k,n}   }
\end{equation*}
expands as
\begin{equation}\label{Expand}
	\frac{  \frac{X_{n,n}- U(n/k)}{a_0(n/k)}  -\frac{1}{\log 2}P_n - \frac{1}{\log 2} \intab{1/2}{1} \frac{U\bigl(n/(2ks) \bigr)- U(n/k)}{a_0(n/k)}\, \frac{ds}{s} }{  Q_n  -  \frac{X_{n-2k,n}- U(n/k)}{a_0(n/k)}},
\end{equation}
where $P_n$ and $Q_n$ are defined and accounted for in Lemma \ref{LemVector}. Under the stated conditions in the Theorem, in particular condition \eqref{2ERVU} of regular variation of second order,  we have for the remainder building blocks:
\begin{eqnarray*}
	\frac{X_{n,n}- U(n/k)}{a_0(n/k)} & \id &\left\{ \begin{array}{ll}
                                \frac{k^{\gamma}(Y_{n,n}/n)^{\gamma}-1}{\gamma}\bigl(1+o_p(1)\bigr), & \mbox{ }  \gamma \neq 0, \\
                              \bigl(\log(Y_{n,n}/n) + \log k\bigr)\bigl(1+o_p(1)\bigr) , & \mbox{ }  \gamma= 0,
			 \end{array}
                            \right.\\
	\frac{X_{n-2k,n}- U(n/k)}{a_0(n/k)} & \id &\left\{ \begin{array}{ll}
                                \frac{2^{-\gamma}-1}{\gamma} + O_p\bigl( 1/\sqrt{k}\bigr), & \mbox{ }  \gamma \neq 0, \\
                            \log (1/2) + O_p\bigl( 1/\sqrt{k}\bigr), & \mbox{ }  \gamma= 0
			 \end{array}
                            \right.
\end{eqnarray*}
(cf. proof of Proposition \ref{PropxF}), and
\begin{equation*}
	b(\gamma):= - \frac{1}{\log 2} \intab{1/2}{1} \frac{U\bigl(n/(2ks) \bigr)- U(n/k)}{a_0(n/k)}\, \frac{ds}{s} = \left\{ \begin{array}{ll}
                               \frac{\gamma \log 2 -1 +2^{-\gamma}}{\gamma^2 \log 2} + O\bigl( A_0(n/k)\bigr), & \mbox{ }  \gamma \neq 0, \\
                               \frac{\log 2}{2} + O\bigl( A_0(n/k)\bigr), & \mbox{ }  \gamma= 0.
			 \end{array}
                            \right.
\end{equation*}
Plugging all the blocks above back in expression $\eqref{Expand}$ for $G_{n,k}$, we therefore obtain:

\noindent if $\gamma=0$,
\begin{eqnarray*}
	G^*_{n,k}(0)&=& \log 2 \,G_{n,k} - \bigl( \log k +\frac{\log 2}{2}\bigr) \\
	&=& \log 2 \frac{ \log (Y_{n,n}/n) +b(0)-\frac{\log 2}{2} +O_p\bigl( \frac{\log k}{\sqrt{k}}\bigr)+ O_p\bigl( A_0(n/k)\bigr)}{\log 2 +  O_p\bigl( \frac{1}{\sqrt{k}}\bigr)},
\end{eqnarray*}
whereas, if $\gamma \neq 0$,
\begin{eqnarray*}
	G^*_{n,k}(0)&=& \log 2 \,G_{n,k} - \bigl( \log k +\frac{\log 2}{2}\bigr) \\
	&=& \log 2 \frac{\frac{k^{\gamma}(Y_{n,n}/n)^{\gamma}-1}{\gamma} +b(\gamma)-\frac{\log 2}{2} -\log k+O_p\bigl( \frac{\log k}{\sqrt{k}}\bigr)+ O_p\bigl( A_0(n/k)\bigr)}{  \frac{2^{-\gamma}-1}{\gamma}  +  O_p\bigl( \frac{1}{\sqrt{k}}\bigr)}.
\end{eqnarray*}
Finally, since $Y_{n,n}/n$ is a non-degenerate random variable, eventually, as it converges to a unit Fr\'echet, the statement follows for $\gamma\in \real$.
\qed

%\newpage
 \section{APPENDIX B: Finite sample properties of $\hat x^F_{RB2}$, $-1/2<\gamma<0$}
 \label{RB2_MC}
The illustrative example about the supercentenarian women's records in section \ref{CaseStudy} also shows how the second order expansion of the general endpoint estimator can be used to remove some contribution to the asymptotic bias. This is the idea underpinning the reduced bias version \eqref{RB2}, provided the true negative EVI stays above  $-1/2$.

 \begin{figure}
    \centering
       \includegraphics[width=0.45\textwidth]{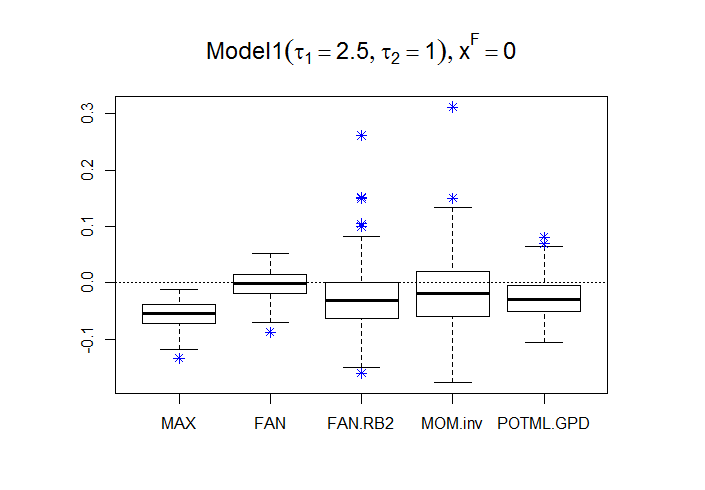}\hfil
     \includegraphics[width=0.45\textwidth]{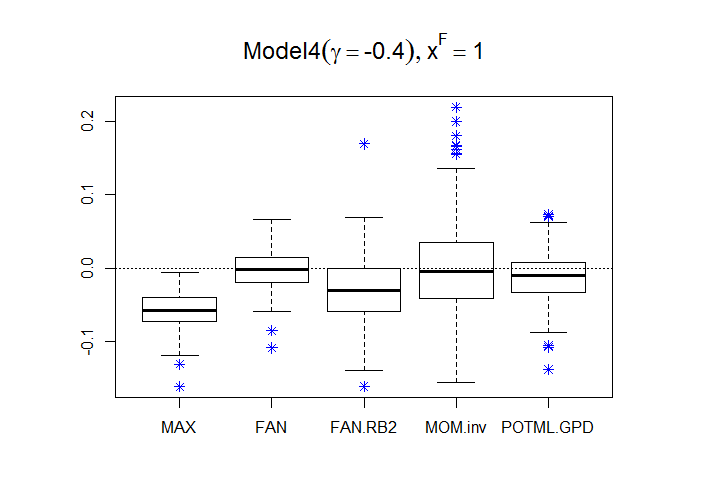}\\
        \includegraphics[width=0.45\textwidth]{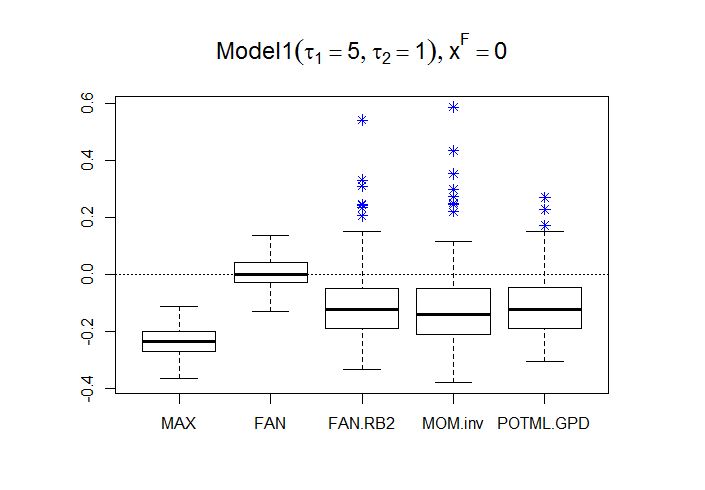}\hfil
        \includegraphics[width=0.45\textwidth]{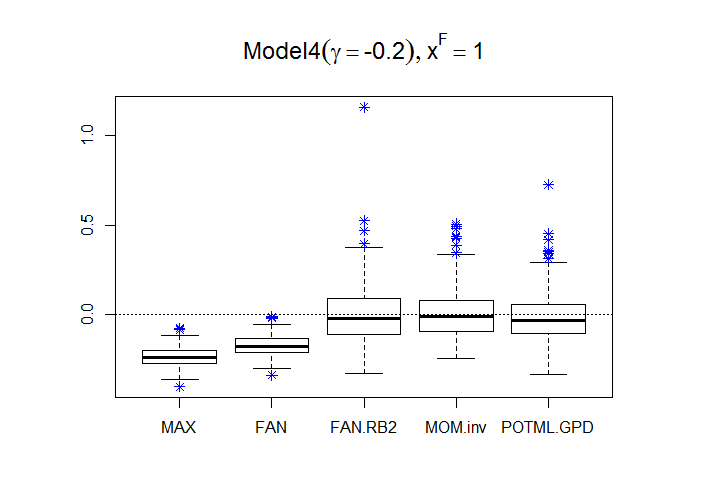}\\
        \includegraphics[width=0.45\textwidth]{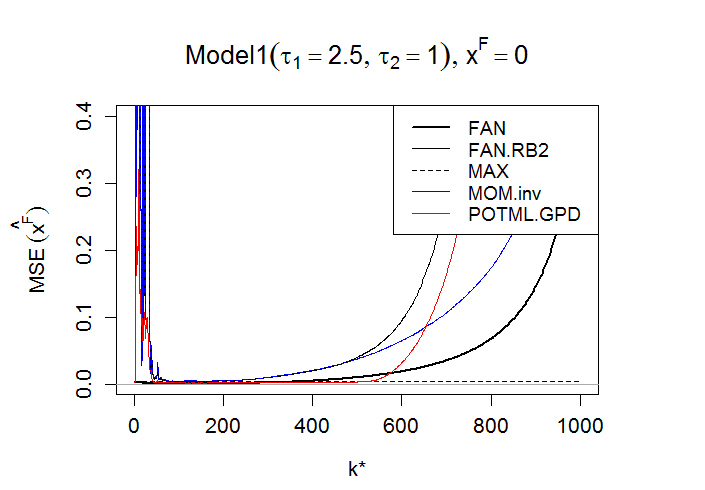}\hfil
     \includegraphics[width=0.45\textwidth]{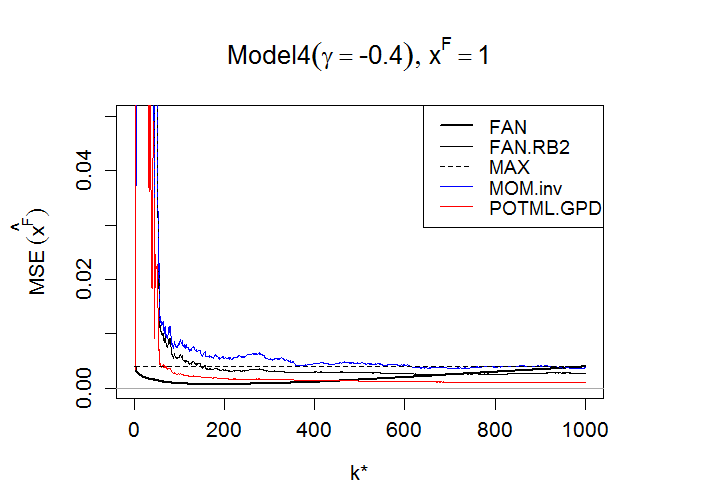}\\
        \includegraphics[width=0.45\textwidth]{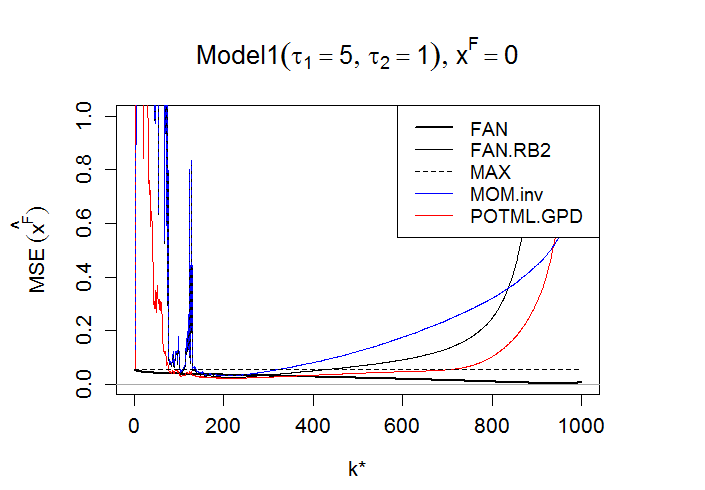}\hfil
        \includegraphics[width=0.45\textwidth]{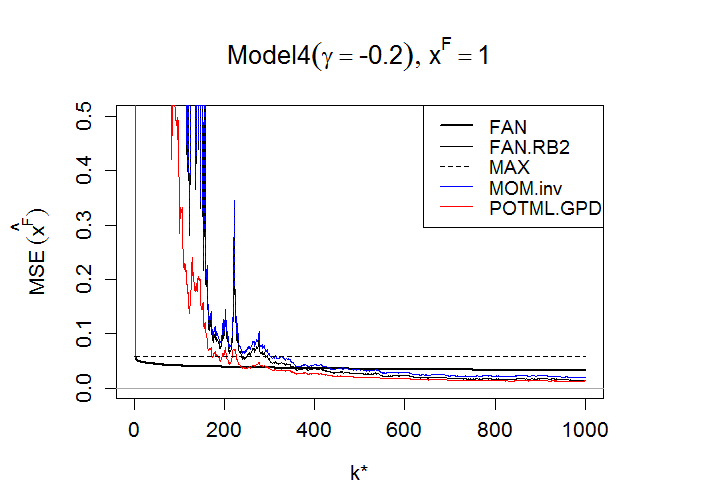}\\
   \caption{\footnotesize Boxplots of the errors $\varepsilon(j,k^*_0),\,\,j=1,\ldots, N=300$ (\emph{top}) and Mean Squared Errors (MSE) plotted against $k^*$, $k^*\leq n$, (\emph{down}) for  MAX, FAN, FAN.RB2, MOM.inv and POTML.GPD endpoint estimators.}
      \label{RB2_box_E}
\end{figure}
Some finite sample results for  $\hat x^F_{RB2}$ are displayed in  Figure \ref{RB2_box_E}. We have generate $N=300$ samples of size $n=1000$ from the parent Models 1 ($x^{F_1}=0$) and 4 ($x^{F_4}=1$). These models were introduced in the simulation study comprising section \ref{SecSims}. Here, we have chosen to set the EVI at the values $-0.4$ and $-0.2$. The practical application in section \ref{CaseStudy} allows to foresee (cf. Figure \ref{xF_data1}) that by reducing the bias in the general endpoint estimator, we end up with a new estimator with larger variance. To this extent, we have furthermore anticipated a new estimator with very similar features to the designated MOM.inv and POTML.GPD endpoint estimators. Now, the simulation results seem to support our ``guess''. The comparative box-plots in Figure \ref{RB2_box_E} show a close resemblance of patterns within the group encompassing the three estimators FAN.RB2, MOM.inv and POTML.GPD,  although there are situations in which the reduced bias version can serve as a good complement to the MOM.inv and POTML.GPD,  particularly for the cases of anomalous behaviour of the likelihood surface, often encountered for the GPD.
 Figure \ref{RB2_box_E} also illustrates the distinctive behavior of the general endpoint estimator, emphasizing lower MSE delivered by this estimator.

%============================================

%==================================================================
%\newpage
 %%%%%%%%%%%%%%%%%%%%%%%%%%%%%%%%%%%%%%%%%%%%%%%%%%%% %%%%%%

\bibliography{bibEndpointNeg}
 \bibliographystyle{apalike}
%\begin{thebibliography}{}
% etc
%\end{thebibliography}

\end{document}